\nonstopmode
\documentclass[10pt]{amsart}
\usepackage{latexsym}
\usepackage{fancyhdr}
\usepackage{amsmath, amssymb}
\usepackage[ansinew]{inputenc}
\usepackage[all]{xy}
\usepackage{pdflscape}
\usepackage{longtable}
\usepackage{rotating}
\usepackage{verbatim}
\usepackage{hyperref}
\usepackage{subfigure}
\usepackage{mathrsfs}
\swapnumbers
\theoremstyle{plain}
\newtheorem*{lemma*}{Lemma}
\newtheorem{lemma}[subsection]{Lemma}
\newtheorem*{theorem*}{Theorem}
\newtheorem{theorem}[subsection]{Theorem}
\newtheorem*{proposition*}{Proposition}
\newtheorem{proposition}[subsection]{Proposition}
\newtheorem*{corollary*}{Corollary}
\newtheorem{corollary}[subsection]{Corollary}
\newtheorem*{claim*}{Claim}
\newtheorem{claim}[subsection]{Claim}

\newtheorem*{conjecture*}{Conjecture}
\newtheorem{question}[subsection]{Question}
\newtheorem*{question*}{Question}
\theoremstyle{definition}
\newtheorem*{definition*}{Definition}

\newtheorem*{example*}{Example}
\newtheorem{example}[subsection]{Example}
\newtheorem{examples}[subsection]{Examples}

\newtheorem*{algorithm*}{Algorithm}
\newtheorem*{remark*}{Remark}
\newtheorem*{remarks*}{Remarks}
\newtheorem{remark}[subsection]{Remark}

\newtheorem*{convention*}{Convention}

\makeatletter
  \let\c@equation\c@subsection
  
\makeatother

\newenvironment{demo}[1]{\par\smallskip\noindent{\bf #1.}}{\par\smallskip}
\numberwithin{equation}{section}

\def\al{\alpha}
\def\be{\beta}
\def\ga{\gamma}
\def\de{\delta}
\def\ep{\epsilon}

\def\ze{\zeta}

\def\th{\theta}

\def\la{\lambda}

\def\rh{\rho}

\def\si{\sigma}

\def\vh{\varphi}
\def\ch{\chi}
\def\ps{\psi}
\def\om{\omega}

\def\De{\Delta}

\def\La{\Lambda}

\def\Ph{\Phi}

\def\C{\mathbb{C}}

\def\I{\mathbb{I}}
\def\K{\mathbb{K}}
\def\N{\mathbb{N}}

\def\R{\mathbb{R}}

\def\cC{\mathcal{C}}

\def\cF{\mathcal{F}}

\def\cH{\mathcal{H}}

\def\cL{\mathcal{L}}

\def\fa{\mathfrak{a}}
\def\fb{\mathfrak{b}}

\def\fs{\mathfrak{s}}

\def\sC{\mathscr{C}}

\def\sE{\mathscr{E}}

\def\p{\partial}

\renewcommand{\Re}{\mathrm{Re}}
\renewcommand{\Im}{\mathrm{Im}}
\def\<{\langle}
\def\>{\rangle}
\renewcommand{\o}{\circ}

\def\CQ{{\cC_Q}}

\let\on=\operatorname
\newcommand{\sr}[1]%
{\ifmmode{}^\dagger\else${}^\dagger$\fi\ifvmode
\vbox to 0pt{\vss
 \hbox to 0pt{\hskip\hsize\hskip1em
 \vbox{\hsize3cm\raggedright\pretolerance10000
 \noindent #1\hfill}\hss}\vss}\else
 \vadjust{\vbox to0pt{\vss%
 \hbox to 0pt{\hskip\hsize\hskip1em%
 \vbox{\hsize3cm\raggedright\pretolerance10000%
 \noindent #1\hfill}\hss}\vss}}\fi%
}

\providecommand{\mapsfrom}{\kern.2em%
\setbox0=\hbox{$\leftarrow$\kern-.10em\rule[0.26mm]{0.1mm}{1.3mm}}\box0%
\kern.3em}
\def\i{^{-1}}

\title[]
{Perturbation theory for normal operators}

\author[A.~Rainer]{Armin Rainer}

\address{A.~Rainer: Fakult\"at f\"ur Mathematik, Universit\"at Wien, 
Nordbergstrasse~15, A-1090 Wien, Austria}
\email{armin.rainer@univie.ac.at}

\begin{document}

\begin{abstract} 
  Let $E \ni x\mapsto A(x)$ be a $\sC$-mapping with values unbounded normal
  operators with common domain of definition and compact resolvent. 
  Here $\sC$ stands for $C^\infty$, $C^\om$ (real analytic), $C^{[M]}$ (Denjoy--Carleman of Beurling or Roumieu type), 
  $C^{0,1}$ (locally Lipschitz), or $C^{k,\al}$.
  The parameter domain $E$ is either $\mathbb R$ or $\mathbb R^n$ or an infinite dimensional convenient vector space.
  We completely describe the $\sC$-dependence on $x$ of the eigenvalues and the eigenvectors of $A(x)$.
  Thereby we extend previously known results for self-adjoint operators to normal operators, partly improve them, 
  and show that they are best possible. 
  For normal matrices $A(x)$ we obtain partly stronger results.
\end{abstract}

\thanks{Supported by the Austrian Science Fund (FWF), Grant P~22218-N13}
\keywords{Perturbation theory, differentiable and Lipschitz eigenvalues and eigenvectors, 
normal unbounded operators, resolution of singularities, Denjoy--Carleman classes}
\subjclass[2010]{26C10, 26E10, 32B20, 47A55}
\date{April 2, 2012}

\maketitle

\section{Introduction and main results} \label{sec:intro}

The purpose of this paper is to prove the following theorem.

\begin{theorem} \label{main}
  Let $x\mapsto A(x)$ be a parameterized family of unbounded normal 
  operators in a Hilbert space $H$ with common domain of definition and with 
  compact resolvent.
  \begin{enumerate} 
  \item[\thetag{A}\phantomsection\label{A}] If $A(x)$ is $C^\infty$ (resp.\ $C^{[M]}$) in $x\in \R$ and 
        if the order of contact of any two unequal eigenvalues is finite at each $x \in \R$,
        then the eigenvalues and the eigenvectors of $A(x)$ admit global $C^\infty$ (resp. $C^{[M]}$) parameterizations in $x$.
        The latter condition is trivially satisfied if $C^{[M]}$ is quasianalytic.  
  \item[\thetag{B}\phantomsection\label{B}] Assume that $C^{[M]}$ is quasianalytic. If $A(x)$ is $C^{[M]}$ in $x\in \R^n$, then
        for each $x_0\in \R^n$ and for each eigenvalue $z$ 
        of $A(x_0)$, there exist a neighborhood $D$ of $z$ in $\C$, a neighborhood $W$ of $x_0$ in $\R^n$, and
        a finite covering $\{\pi_k : U_k \to W\}$ of $W$ by composites of finitely many local blow-ups,
        such that the eigenvalues of $A(\pi_k(y))$ in $D$ and the corresponding eigenvectors can be chosen $C^{[M]}$ in $y$.
  \item[\thetag{C}\phantomsection\label{C}] Assume that $C^{[M]}$ is quasianalytic. If $A(x)$ is $C^{[M]}$ in $x\in \mathbb R^n$, 
        then for each $x_0\in \R^n$ and for each eigenvalue $z$ 
        of $A(x_0)$, there exists a neighborhood $D$ of $z$ in $\C$, such that
  			the eigenvalues of $A(x)$ in $D$ can be parameterized by functions 
        which are locally `piecewise Lipschitz continuous', i.e., belong to $\cL^{C^{[M]}}_{\on{loc}}$ (cf.\ \ref{classL}).
        In particular, they are $SBV_{\on{loc}}$-functions whose classical gradient exists almost everywhere and is locally bounded.  
  \item[\thetag{D}\phantomsection\label{D}] If $x \mapsto A(x)$ is $C^{0,1}$ in $x \in E$, where $E$ is a convenient vector space, 
        then each continuous eigenvalue $E \supseteq U \ni x \mapsto \la(x)$, 
        for $c^\infty$-open $U \subseteq E$, of $A(x)$ is $C^{0,1}$ in $x$. 
        If $x_0 \in E \cap \overline U$ and $c : \R \to E$ is a $C^\infty$-curve with $c(0)=x_0$ and $c((0,1]) \subseteq U$, 
        then $\la \o c|_{(0,1]}$ is globally Lipschitz on $(0,1]$.
        If $E=\R$, then the eigenvalues admit a $C^{0,1}$-parameterization in $x$.    
  \item[\thetag{E}\phantomsection\label{E}] If $x \mapsto A(x)$ is $C^{1,\al}$ in $x \in \R$, for some $\al>0$, 
        then the eigenvalues admit a $C^1$-parameterization in $x$.      
  \item[\thetag{F}\phantomsection\label{F}] If $x \mapsto A(x)$ is $C^{2,\al}$ in $x \in \R$, for some $\al>0$, 
        then the eigenvalues admit a twice differentiable parameterization in $x$.
  \end{enumerate}
\end{theorem}

Let us define the involved notions and explain the results.

\subsection{Definitions and remarks} \label{ssec:def}
For a sequence $M=(M_k)_{k \in \N}$ of positive real numbers, $U \subseteq \R^n$ open, $K \subseteq U$ compact, and $\rh>0$,
consider the set 
\begin{equation} \label{DCest}
  \Big\{\frac{\p^\al f(x)}{\rh^{|\al|} \, |\al|! \, M_{|\al|}} : x \in K, \al \in \N^n \Big\},
\end{equation}
and define the \emph{Denjoy--Carleman classes} 
\begin{align*}
  C^{(M)}(U) &:= \{f \in C^\infty(U) : \forall \text{ compact } K \subseteq U ~\forall \rh>0: \eqref{DCest} \text{ is bounded} \}, \\
  C^{\{M\}}(U) &:= \{f \in C^\infty(U) : \forall \text{ compact } K \subseteq U ~\exists \rh>0: \eqref{DCest} \text{ is bounded} \}.
\end{align*}
The elements of $C^{(M)}(U)$ are said to be of \emph{Beurling type}; those of $C^{\{M\}}(U)$ of \emph{Roumieu type}. 
If $M_k=1$, for all $k$, then $C^{(M)}(U)$ consists of the restrictions to $U$ of the real and imaginary parts of all entire functions, 
while $C^{\{M\}}(U)$ coincides with the ring $C^\om(U)$ of real analytic functions on $U$. 

We use the notation $C^{[M]}$ for either $C^{(M)}$ or $C^{\{M\}}$ with the following restriction: 
Statements that involve more than one $C^{[M]}$ symbol must not be interpreted by mixing $C^{(M)}$ and $C^{\{M\}}$.

We shall always assume that $M=(M_k)$ has the following regularity properties:
\begin{enumerate}
  \item[\thetag{M$_1$}\phantomsection\label{M_1}] Log-convexity: $M_k^2 \le M_{k-1} \, M_{k+1}$ for all $k$.
  \item[\thetag{M$_2$}\phantomsection\label{M_2}] Stability under derivation: $\sup_{k} \big(\frac{M_{k+1}}{M_k}\big)^{1/k} < \infty$.
\end{enumerate}
Then $C^{[M]}$ is stable under composition and derivation.   
Moreover, $C^{\{M\}}\supseteq C^\om$, the $C^{\{M\}}$ inverse function theorem holds, and $C^{\{M\}}$ is closed under solving ODEs. 
The $C^{(M)}$ inverse function theorem is valid and $C^{(M)}$ is closed under solving ODEs if additionally 
$M_{k+1}/M_k \to \infty$. This is satisfied if
\begin{enumerate}
  \item[\thetag{M$_3$}\phantomsection\label{M_3}] $M_k^{1/k} \to \infty$
\end{enumerate}
which will always be assumed in the Beurling case. 
Condition \thetag{\hyperref[M_3]{M$_3$}} is equivalent to $C^\om \subseteq C^{(M)}$ and in turn to $C^\om \subsetneq C^{\{M\}}$.
The classes $C^{[M]}$ are quasianalytic, i.e., infinite Taylor expansion is injective, if and only if the following condition holds:
\begin{enumerate}
  \item[\thetag{M$_4$}\phantomsection\label{M_4}] Quasianalyticity: $\sum_{k} \frac{M_k}{(k+1)M_{k+1}}=\infty$. 
\end{enumerate}
For more details on Denjoy--Carleman classes see \cite{Thilliez08}, \cite{KMRc}, \cite{KMRq}, \cite{KMRu}, and references therein.

A \emph{convenient} vector space is a real locally convex vector space
$E$ satisfying the following equivalent conditions: Mackey Cauchy sequences
converge; $C^\infty$-curves in $E$ are locally integrable in $E$; 
a curve $c:\mathbb R\to E$ is $C^\infty$ if and only if 
$\ell\circ c$ is $C^\infty$
for all continuous linear functionals $\ell$.
The \emph{$c^\infty$-topology} on $E$
is the final topology with respect to all $C^\infty$-curves.
Functions $f$ defined on $c^\infty$-open 
subsets of convenient vector spaces $E$ 
are called $C^{k,\al}$ if $f\circ c$ is $C^{k,\al}$ for every $C^\infty$-curve
$c$. If $E$ is a Banach space, then a $C^{k,\al}$-function is $C^k$ and the $k$th derivative is locally H\"older continuous
of order $\al$ in the usual sense. This has been proved in \cite{Faure91}, see also the lemma in \cite{KMR}. 
For the Lipschitz case see \cite{FK88} and \cite[12.7 and 12.8]{KM97}.

That $A(x)$ is a $C^\infty$, $C^{[M]}$, or a $C^{k,\al}$-family of 
unbounded normal operators means the following:
There is a dense subspace $V$ of the Hilbert space $H$ 
so that $V$ is the domain of definition of each $A(x)$, $A(x)$ has closed graph, and we have
$A(x)A(x)^*=A(x)^*A(x)$ wherever defined. 
Moreover, we require 
that $x\mapsto \langle A(x)u \mid v\rangle$ 
is $C^\infty$, $C^{[M]}$, or $C^{k,\al}$, for each $u\in V$ and $v\in H$. 
This implies that $x\mapsto A(x)u$ is of the same class as a mapping $E \to H$ (where $E$ is either $\mathbb R$ 
or $\mathbb R^n$ or an infinite dimensional convenient vector space)
for each $u\in V$, 
by \cite[2.3]{KM97} for $C^\infty$, 
by \cite[4.3, 4.4, 4.5, and 5.1]{KMRu} 
for $C^{[M]}$,
and by \cite[2.3]{KM97}, \cite[2.6.2]{FK88} or \cite[4.14.4]{Faure91} 
for $C^{k,\al}$, because $C^{k,\al}$ can be described by boundedness conditions only 
and for these the uniform boundedness principle is valid.
Note that the real analytic case is included since $C^\om=C^{\{(1)_k\}}$.

If $A$ depends on a single real parameter $x$, then the eigenvalues of $A$ may be chosen continuously
near each $(x_0,z)$, where $z$ is an eigenvalue of $A(x_0)$, 
see \cite[II Thm.\ 5.2]{Kato76}. The order of vanishing of a continuous function germ $f$ at $0 \in \R$  
is the supremum of all integers $p$ such that $f(x)=x^p g(x)$, where $g$ is continuous; likewise at any $x_0 \in \R$.
The order of contact of two continuous function germs   
is the order of vanishing of their difference.  

A \emph{local blow-up $\Ph$} over an open subset $U$ of a $C^{[M]}$-manifold $X$ means the composite
$\Ph = \iota \o \vh$ of a blow-up $\vh : U' \to U$ with center a $C^{[M]}$-submanifold and 
of the inclusion $\iota : U \to X$.

A sequence of functions 
$\la_i$ is said to parameterize the eigenvalues of $A$, if,
for each $z\in \mathbb C$, 
the cardinality $|\{i: \la_i(x)=z\}|$ equals
the multiplicity of $z$ as an eigenvalue of $A(x)$. 

An $SBV$-function is a special function of bounded variation, i.e., a function having bounded variation whose distributional derivative 
has trivial Cantor part, see \cite{AmbrosioDeGiorgi88} and \cite{AFP00}.

\subsection{Explanation of the results and background}
The novelty of the results in Theorem~\ref{main} and of the partly stronger 
finite dimensional versions of \thetag{\hyperref[A]{A}}--\thetag{\hyperref[F]{F}} for normal matrices which will be shown in 
the course of the proof of Theorem~\ref{main}    
is threefold:
\begin{itemize}
  \item The results are well-known if all operators $A(x)$ are self-adjoint; at least in some weaker formulation.
  We show that the assumption of self-adjointness can be replaced by normality, essentially without changing the conclusions 
  (only in \thetag{\hyperref[D]{D}} we additionally have to assume continuity if $\dim E > 1$).
  \item We achieve utmost generality, at least for matrices, by working in abstractly defined quasianalytic subclasses of 
  $C^\infty$ which present a minimal setting for our method of proof. For unbounded operators we restrict to $C^{[M]}$. 
  \item We partly even improve the results for self-adjoint operators and show that they are then best possible. 
\end{itemize}

Let us briefly describe what was previously known.
If all operators $A(x)$ are self-adjoint, then \thetag{\hyperref[A]{A}} is due to Rellich \cite{Rellich42V} in 1942 for $C^\om$, 
to \cite{AKLM98} for $C^\infty$, 
and to \cite{KMRp} for $C^{\{M\}}$ (with special $M=(M_k)$);
the normal case follows for $C^\om$ from an observation due to Butler, see \cite[II Thm.\ 1.10]{Kato76} and \cite[3.5.1]{Baumgaertel85}.  
Part \thetag{\hyperref[B]{B}} is due to \cite{KurdykaPaunescu08} for $C^\om$-families of symmetric matrices 
and to \cite{KMRp} for unbounded self-adjoint operators; 
in \cite{KMRp} (see also \cite{RainerAC} and \cite{RainerQA}) the normal case is treated, 
but there in addition we had to use local power substitutions.
In the self-adjoint case,
part \thetag{\hyperref[C]{C}} and part \thetag{\hyperref[D]{D}} are consequences of \cite[9.6]{RainerQA} and \cite{KMR}, and 
part \thetag{\hyperref[E]{E}} was proved in \cite{KM03}. 
Part \thetag{\hyperref[F]{F}} was shown in \cite{KM03} under the assumption that $\R \ni x \to A(x)$ is a $C^\infty$-curve 
(or, more precisely, $C^{3n,\al}$, if the multiplicity of an eigenvalues does never exceed $n$) 
of self-adjoint operators. Our proof of \thetag{\hyperref[F]{F}} works for normal $A$ and needs only the assumption $C^{2,\al}$.

It is somewhat surprising that these results carry over to normal operators. 
For Hermitian matrices the characteristic polynomial is hyperbolic, i.e., all its roots are real, and 
the roots of families of hyperbolic polynomials admit `nice' parameterizations, which are reflected by the regularity 
properties of the eigenvalues and the eigenvectors. 
For instance, the roots of a hyperbolic polynomial with coefficients in some quasianalytic class of functions admit parameterizations
in the same class after desingularization by means of local blow-ups (of the parameter space), 
see \cite{RainerQA} and \cite{KurdykaPaunescu08} for $C^\om$; 
and the (increasingly ordered) roots are locally Lipschitz, provided
that the coefficients are in $C^n$, where $n$ is the degree, see \cite{Bronshtein79}. 
The perturbation theory for complex polynomials is considerably weaker: In general, 
local power substitutions are needed in order to desingularize, and the roots cannot satisfy a 
local Lipschitz condition, e.g., $z^2-x =0$, $x \in \R$, see \cite{RainerQA}.
However, not every quasianalytic family of polynomials appears as the characteristic polynomial of a quasianalytic family of normal matrices.
In fact, the set of normal complex $n \times n$ matrices forms a real $n^2+n$ dimensional stratified submanifold of $\R^{2n^2}$ 
(the set of all complex $n \times n$ matrices), see e.g.\ \cite{Huhtanen01}.  
So the normality condition implies perturbation results for operators stronger than predicted by the perturbation theory for polynomials.

The results in \thetag{\hyperref[B]{B}} and \thetag{\hyperref[C]{C}} seem to be new even in the real analytic setting. 
However, we shall work in a minimal setting making the proofs (in particular desingularization) work, namely   
subclasses of $C^\infty$ which are quasianalytic and have certain stability properties, see Section~\ref{secC}. 
Only when passing to infinite dimensions we will restrict to the framework of Denjoy--Carleman classes 
for which we have developed the required principles of calculus beyond Banach spaces in \cite{KMRc,KMRq,KMRu}.
One may expect analogous results for any suitable quasianalytic function class.

In \thetag{\hyperref[D]{D}} we need to assume continuity of $x \mapsto \la(x)$ if $\dim E > 1$, 
since in general there will not exist continuous parameterizations of the single
eigenvalues, see Example~\ref{excont}. However, it might be that the supplement in \thetag{\hyperref[C]{C}} is still true without that 
assumption,
i.e., that a $C^{0,1}$-family $\R^n \ni x \mapsto A(x)$ of normal complex matrices
admits a parameterization of its eigenvalues by $SBV_{\on{loc}}$-functions whose classical gradient exists a.e.\ and is locally bounded, 
see Question~\ref{question}. 

The conclusions in \thetag{\hyperref[E]{E}} and \thetag{\hyperref[F]{F}} are optimal in the following sense:
There exist $C^\infty$-curves (even non-quasianalytic $C^{[M]}$) of real symmetric $2 \times 2$ matrices whose eigenvalues 
do not admit a parameterization in $C^{1,\al}$ for any $\al>0$, see the examples in \cite{KM03} and \cite{KMRp}.  
We also want to stress that \thetag{\hyperref[A]{A}}, \thetag{\hyperref[E]{E}}, and \thetag{\hyperref[F]{F}} are no longer 
true if the parameter domain has more than one dimension: 
The eigenvalues $\pm \sqrt{x^2+y^2}$ of the real analytic family
$\big(\begin{smallmatrix}
  x & y \\ y & -x
\end{smallmatrix}\big)$, $x,y \in \R$, are not $C^1$ at the origin, see Example~\ref{ex1}.

We point out that the assumptions in Theorem~\ref{main} may be slightly relaxed, 
if all $A(x)$ are m-sectorial operators. 
In that case it suffices to assume that the associated quadratic forms $\fa(x)$ 
have common domain of definition $V$ and $x \mapsto \fa(x)(u)$ is of the respective class for each $u \in V$, see Remark~\ref{m-sec}.

The paper is organized as follows: 
We introduce and describe the classes of smooth functions we shall be working with in Section~\ref{secC} and polynomials with 
coefficients in these classes 
in Section~\ref{secCpoly}.
In Section~\ref{secpoly} we show that a quasianalytic polynomial is solvable (i.e., admits roots in the same class as the coefficients) 
along quasianalytic arcs if and only if it is solvable after blowing up (the parameter space). This will be used in the proof of \thetag{\hyperref[B]{B}}.
We shall prove (partly stronger) finite dimensional versions of \thetag{\hyperref[A]{A}}--\thetag{\hyperref[F]{F}} for normal matrices in 
Section~\ref{ramatrix} and \ref{Lipmatrix}. 
The proof of the Theorem~\ref{main} will finally be completed in Section~\ref{unbounded}. 
Several examples in Section~\ref{optimal} will show that the results are best possible 
in the sense that, generally, the assumptions cannot be weakened and the conclusions cannot be strengthened.
In particular, the results are no longer true if $A$ is a family of merely diagonalizable matrices.

\subsection*{Notation}
The notation $C^{[M]}$ stands for either $C^{(M)}$ or $C^{\{M\}}$ with the following restriction: 
Statements that involve more than one $C^{[M]}$ symbol must not be interpreted by mixing $C^{(M)}$ and $C^{\{M\}}$.

Let $\N = \N_{>0} \cup \{0\}$. 
For $\al=(\al_1,\ldots,\al_q) \in \N^q$ and $x = (x_1,\ldots,x_q) \in \R^q$
we write $\al!=\al_1! \cdots \al_q!$, $|\al|= \al_1 +\cdots+ \al_q$, $x^\al = x_1^{\al_1}\cdots x_q^{\al_q}$, and 
$\p^\al=\p^{|\al|}/\p x_1^{\al_1} \cdots \p x_q^{\al_q}$. 
We shall also use $\p_i = \p/\p x_i$, $d$ for the Fr\'echet derivative, and $d_v$ for the directional derivative in direction $v$. 
If $\al,\be \in \N^q$, then $\al \le \be$ means $\al_i \le \be_i$ for all $i$.

For a $C^\infty$ function germ $f$ at $a \in \R^q$ we denote by 
$\widehat{f}_a \in \cF_q$ its Taylor series at $a$,
where $\cF_q$ is the ring of formal power series in $q$ variables. 
We write $\cF_q^\K = \K[[x_1,\ldots,x_q]]$ if we want to stress that the coefficients belong to $\K$ 
(where $\K=\R$ or $\K=\C$) and the variables are $x_1,\ldots,x_q$.
We also use $\widehat{f}=\widehat{f}_0$.
We write $\om(F)$ for the order of $F \in \cF_q$, i.e., the lowest degree of non-zero monomials in $F$, 
with the convention $\om(0) = + \infty$. 
For a $C^\infty$ function germ $f$ at $0$ we set $\om(f):=\om(\widehat{f}\,)$. 

$\on{S}_n$ denotes the symmetric group on $\{1,2,\ldots,n\}$. It acts on $\C^n$ by permuting the coordinates: 
$\si.z = (z_{\si(1)},\ldots,z_{\si(n)})$ for $z=(z_1,\ldots,z_n) \in \C^n$ and $\si \in \on{S}_n$. 
This action is denoted by $\on{S}_n : \C^n$. 
The isotropy subgroup that fixes $z$ is denoted by 
$(\on{S}_n)_z = \{\si \in \on{S}_n : \si.z=z\}$.
The elementary symmetric functions $\si_j = \sum_{i_1 <\cdots<i_j} z_{i_1}\cdots z_{i_j}$ generate the algebra of 
symmetric polynomials $\C[\C^n]^{\on{S}_n}$.

We write $|S|$ for the cardinality of a finite set $S$
and denote by $\cH^q$  
the $q$-dimensional Hausdorff  
measure. 

$L(E,F)$ is the space of bounded linear mappings $E \to F$.

\section{Smooth function classes} \label{secC}

\subsection{Classes of \texorpdfstring{$C^\infty$}{Cinfty}-functions} \label{cC}
Let us assume that for every open $U \subseteq \R^q$, $q \in \N$, we have a subalgebra $\cC(U)$ 
of $C^\infty(U)=C^\infty(U,\R)$ so that the following assumptions \thetag{\hyperref[C_1]{$\cC_1$}}--\thetag{\hyperref[C_5]{$\cC_5$}} are satisfied.
\begin{enumerate}
\item[\thetag{$\cC_1$}\phantomsection\label{C_1}]  {\it $\cC$ contains the restrictions of polynomial functions}.
  The algebra of restrictions to $U$ of polynomial functions on $\R^q$ is contained in $\cC(U)$.
\item[\thetag{$\cC_2$}\phantomsection\label{C_2}] {\it $\cC$ is closed under composition}. 
If $V \subseteq \R^p$ is open and $\vh=(\vh_1,\ldots,\vh_p) : U \to V$ is a mapping 
with each $\vh_i \in \cC(U)$, then $f \o \vh \in \cC(U)$, for all $f \in \cC(V)$.
\item[\thetag{$\cC_3$}\phantomsection\label{C_3}] {\it $\cC$ is closed under derivation}. If $f \in \cC(U)$ and $1 \le i \le q$,
then $\p_i f \in \cC(U)$. 
\item[\thetag{$\cC_4$}\phantomsection\label{C_4}] {\it $\cC$ is closed under division by a coordinate}. 
If $f \in \cC(U)$ is identically $0$ along a hyperplane $\{x : x_i=a_i\}$, 
then $f(x)= (x_i-a_i) h(x)$, where $h \in \cC(U)$.
\item[\thetag{$\cC_5$}\phantomsection\label{C_5}] {\it $\cC$ is closed under taking the inverse}. 
Let $\vh : U \to V$ be a $\cC$-mapping between open subsets $U$ and $V$ in $\R^q$. 
Let $a \in U$, $\vh(a)=b$, and suppose that the Jacobian matrix $(\p \vh/\p x)(a)$ is invertible. Then 
there exist neighborhoods $U'$ of $a$, $V'$ of $b$, and a $\cC$-mapping $\ps : V' \to U'$ such that $\ps(b)=a$ and $\vh \o \ps = \on{id}_{V'}$. 
\end{enumerate}
A mapping $\vh : U \to V$ between open subsets $U \subseteq \R^q$ and $V \subseteq \R^p$ 
is called a \emph{$\cC$-mapping} if $f \o \vh \in \cC(U)$, for every $f \in \cC(V)$. 
It follows from \thetag{\hyperref[C_1]{$\cC_1$}} and \thetag{\hyperref[C_2]{$\cC_2$}} that $\vh=(\vh_1,\ldots,\vh_p)$ is a $\cC$-mapping if and only if $\vh_i \in \cC(U)$, 
for all $1 \le i \le p$.

Property \thetag{\hyperref[C_5]{$\cC_5$}} is equivalent to the \emph{implicit function theorem in $\cC$}:
Let $U \subseteq \R^q \times \R^p$ be open. Suppose that $f_1,\ldots,f_p \in \cC(U)$, $(a,b) \in U$, $f(a,b)=0$, 
and $(\p f/\p y)(a,b)$ is invertible, where $f=(f_1,\ldots,f_p)$.
Then there is a neighborhood $V \times W$ of $(a,b)$ in $U$ and a $\cC$-mapping $g : V \to W$ such that $g(a)=b$ and $f(x,g(x))=0$, for $x \in V$.

It follows from \thetag{\hyperref[C_5]{$\cC_5$}} that {\it $\cC$ is closed under taking the reciprocal}: If $f \in \cC(U)$ vanishes nowhere in $U$, 
then $1/f \in \cC(U)$. 

Frequently, we shall also require the following condition.
\begin{enumerate}
\item[\thetag{Q}\phantomsection\label{Q}] {\it $\cC$ is quasianalytic}. If $f \in \cC(U)$ and for $a \in U$ the Taylor series of $f$ at $a$ vanishes (i.e., $\widehat{f}_a=0$) then $f$ vanishes in a neighborhood of $a$. 
\end{enumerate}
Since $\{x : \widehat{f}_x=0\}$ is closed in $U$, condition \thetag{\hyperref[Q]{Q}} is equivalent to the following property: 
If $U$ is connected, then, for each $a \in U$,
the Taylor series homomorphism $\cC(U) \to \cF_q$, 
$f \mapsto \widehat{f}_a$, is injective.

Occasionally, we will need a further condition.
\begin{enumerate}
\item[\thetag{$\cC_6$}\phantomsection\label{C_6}] {\it $\cC$ is closed under solving ODEs}.  
  Let $I\subseteq \R$ be an open interval and let $U \subseteq \R^q$ be open.
  Consider the initial value problem
  \[
  x' = f(t,x), \quad x(0) = y,
  \]
  where $f : I \times U \to \R^q$ is a $\cC$-mapping. Then the smooth solution $x=x(t,y)$ 
  is of class $\cC$ wherever it exists.  
\end{enumerate}

A complex-valued function $f : U \to \C$ is said to be a \emph{$\cC$-function}, or to belong to $\cC(U,\C)$, 
if $(\Re f, \Im f) : U \to \R^2$ is a $\cC$-mapping.
It is immediately verified that \thetag{\hyperref[C_3]{$\cC_3$}} and \thetag{\hyperref[C_4]{$\cC_4$}} also hold for complex-valued functions $f \in \cC(U,\C)$; as well as \thetag{\hyperref[Q]{Q}} if assumed.

\begin{convention*}
  From now on, $\cC$ shall denote a fixed, but arbitrary, class of $C^\infty$-functions satisfying the conditions
  \thetag{\hyperref[C_1]{$\cC_1$}}--\thetag{\hyperref[C_5]{$\cC_5$}}. 
  We shall write $\CQ$ for a class $\cC$ which is required to satisfy \thetag{\hyperref[Q]{Q}}.
  It will be explicitly stated when \thetag{\hyperref[C_6]{$\cC_6$}} is assumed.
\end{convention*}

Note that $\cC$ might be $C^\infty$ and $\CQ$ might be $C^\om$. Here are some more examples.

\begin{examples}[Denjoy--Carleman classes (\cite{Thilliez08}, \cite{KMRc}, \cite{KMRu}, and references therein)] 
  (1) \emph{Denjoy--Carleman classes of Roumieu type}:
  If $M=(M_k)$ is a positive log-convex sequence which is stable under derivation 
  (see \thetag{\hyperref[M_1]{M$_1$}} and \thetag{\hyperref[M_2]{M$_2$}}),
  then the Denjoy--Carleman class of Roumieu type $C^{\{M\}}$ has the properties  
  \thetag{\hyperref[C_1]{$\cC_1$}}--\thetag{\hyperref[C_6]{$\cC_6$}}; 
  see \cite[Section 4]{BM04} for \thetag{\hyperref[C_1]{$\cC_1$}}--\thetag{\hyperref[C_5]{$\cC_5$}} and \cite{Komatsu80} 
  for \thetag{\hyperref[C_6]{$\cC_6$}}. 
  In particular, this is true for all \emph{Gevrey classes} $G^{1+s}=C^{\{(k!^s)_k\}}$, $s\ge 0$.
  If $M=(M_k)$ additionally satisfies \thetag{\hyperref[M_4]{M$_4$}}, then $C^{\{M\}}$ is quasianalytic \thetag{\hyperref[Q]{Q}}.
  Among the Gevrey classes only $G^1 = C^\om$ has this property. 
  However, by setting
  \[
  M_k^{\de,n} := \frac1{k!}\big(k \cdot \log(k) \cdot \cdots \cdot \log^{n-1}(k) \cdot (\log^n(k))^\de\big)^k,
  \] 
  where $\log^n$ denotes the $n$-fold composition of $\log$, we obtain for each $0< \de \le 1$ and each $n \in \N_{>0}$
  a quasianalytic class $C^{\{M^{\de,n}\}}$ satisfying all required conditions, 
  and $C^{\{M^{\de,n}\}} \ne C^{\{M^{\de',n'}\}}$ if $(\de,n) \ne (\de',n')$; 
  see \cite[1.9]{KMRq}.

  (2) \emph{Denjoy--Carleman classes of Beurling type}:
  If $M=(M_k)$ is a positive log-convex sequence which is stable under derivation 
  (see \thetag{\hyperref[M_1]{M$_1$}} and \thetag{\hyperref[M_2]{M$_2$}}),
  then the Denjoy--Carleman class of Beurling type $C^{(M)}$ has the properties  
  \thetag{\hyperref[C_1]{$\cC_1$}}--\thetag{\hyperref[C_4]{$\cC_4$}}. 
  Properties \thetag{\hyperref[C_5]{$\cC_5$}} and \thetag{\hyperref[C_6]{$\cC_6$}} are satisfied 
  if additionally $M_{k+1}/M_k \to \infty$ (which follows from \thetag{\hyperref[M_3]{M$_3$}}). 
  See \cite[2.1]{KMRu} for references.
  Again the non-quasianalytic classes $C^{((k!^s)_k)}$, $s>0$, and the quasianalytic classes $C^{(M^{\de,n})}$ have all required properties.
  
  If $C^{[M]}$ has all properties \thetag{$\cC_i$} but \thetag{\hyperref[C_3]{$\cC_3$}}, i.e., it
  is not closed under derivations, 
  then $\bigcup_{j \in \N} C^{[M^{+j}]}$, where $M^{+j}_k:= M_{k+j}$, has the properties 
  \thetag{\hyperref[C_1]{$\cC_1$}}--\thetag{\hyperref[C_5]{$\cC_6$}}, and, moreover, 
  it satisfies \thetag{\hyperref[Q]{Q}} if and only if $C^{[M]}$ does.
\end{examples}

\subsection{Resolution of singularities in \texorpdfstring{$\CQ$}{C_Q}} \label{Cmf}

A $\cC$-manifold is a $C^\infty$-manifold such that all chart change mappings are of class $\cC$. 
This provides a category $\underline{\cC}$ of \emph{$\cC$-manifolds} and 
\emph{$\cC$-mappings}.

The implicit function property \thetag{\hyperref[C_5]{$\cC_5$}} implies that a \emph{smooth} (i.e., not singular) 
subset of a $\cC$-manifold is a 
$\cC$-submanifold:
Let $M$ be a $\cC$-manifold. Suppose that $U$ is open in $M$, $g_1,\ldots,g_p \in \cC(U)$, and the gradients 
$\nabla g_i$ are linearly independent at every point of the zero set $X:=\{x \in U : g_i(x)=0 \text{ for all }i\}$. 
Then $X$ is a closed $\cC$-submanifold of $U$ of codimension $p$. 

The category $\underline \cC$ 
is closed under blowing up with center a closed 
$\cC$-submanifold. 

We shall use a simple version of the desingularization theorem of Hironaka \cite{Hironaka64} for $\CQ$-function classes due 
to Bierstone and Milman \cite{BM97, BM04}.
We use the terminology therein.

\begin{theorem}[{\cite[5.12]{BM04}}] \label{resth}
  Let $M$ be a $\CQ$-manifold, $X$ a closed $\CQ$-hypersurface in $M$, and $K$ a compact subset of $M$. 
  Then, there is a neighborhood $W$ of $K$ and a surjective mapping $\vh : W' \to W$ of class $\CQ$, such that:
  \begin{enumerate}
    \item $\vh$ is a composite of finitely many $\CQ$-mappings, 
    each of which is either a blow-up with smooth center 
    (that is nowhere dense in the smooth points of the strict transform of $X$) 
    or a surjection of the form $\bigsqcup_j U_j \to \bigcup_j U_j$, 
    where the latter is a finite covering of the target space by coordinate charts.
    \item The final strict transform $X'$ of $X$ is smooth, and $\vh^{-1}(X)$ has only normal crossings. 
    (In fact $\vh^{-1}(X)$ and $\det d \vh$ simultaneously have only normal crossings, 
    where $d \vh$ is the Jacobian matrix of $\vh$ with respect to any local coordinate system.)
  \end{enumerate}
\end{theorem}

See \cite[5.9 and 5.10]{BM04} and \cite{BM97} for stronger desingularization theorems in $\CQ$. 

A real- or complex-valued $\CQ$-function on a $\CQ$-manifold $M$ is said to 
have only \emph{normal crossings} if each point in $M$ admits a coordinate neighborhood $U$ with coordinates 
$x=(x_1,\ldots,x_q)$ such that
\[
f(x)=x^\al g(x), \quad x \in U,
\]
where $g$ is a non-vanishing $\CQ$-function on $U$, and $\al \in \N^q$.
Observe that, if a product of $\CQ$-functions has only normal crossings, then each factor has only normal crossings.

Let $f \in \CQ(M,\C)$ and let $K \subseteq M$ be compact. Then there exists a neighborhood $W$ of $K$ 
and a finite covering $\{\pi_k : U_k \to W\}$ of $W$ by $\CQ$-mappings $\pi_k$,
each of which is a composite of finitely many local blow-ups, such that, for each $k$, 
the function $f \o \pi_k$ has only normal crossings. This follows from Theorem~\ref{resth} applied to the real-valued $\CQ$-function 
$|f|^2=f\overline f$ and from the previous observation.

By a \emph{local blow-up} $\Ph$ over an open subset $U$ of a $\CQ$-manifold $M$ we mean 
the composite
$\Ph = \iota \o \vh$ of a blow-up $\vh : U' \to U$ with smooth center and 
of the inclusion $\iota : U \to M$.

We shall need the following well-known lemma.

\begin{lemma}[{\cite[7.7]{BM04}}, {\cite[4.7]{BM88}}, or {\cite[6.3]{RainerQA}}] \label{order} 
Let $\al,\be,\ga \in \N^q$ and let $a,b,c$ be non-vanishing germs of real- or complex-valued $\CQ$-functions at the origin of 
$\R^q$. If 
$x^\al a(x) - x^\be b(x) = x^\ga c(x)$,
then either $\al \le \be$ or $\be \le \al$. 
\end{lemma}

The following simple observation will be used repeatedly.

\begin{lemma} \label{unique}
  Let $I \subseteq \R$ be an open interval. 
  Let $f_j,g_j : I \to \C$, $1 \le j \le n$, be $\cC$-functions such that 
  $|\{j : f_j(t)=z\}| = |\{j : g_j(t)=z\}|$ for all $t \in I$ and $z \in \C$.
  Assume that at each $t_0 \in I$ the order of contact of any two elements of $\{f_j\}$ (equivalently $\{g_j\}$) is finite unless their germs at $t_0$ coincide. 
  Then $\{f_j\}$ and $\{g_j\}$ differ by a constant permutation.  
\end{lemma}

The assumption on the order of contact is trivially satisfied if the functions are of class $\CQ$.

\begin{demo}{Proof}
  Set $f=(f_1,\ldots,f_n)$, $g=(g_1,\ldots,g_n)$, and consider the set
  \[
  J:= \{t \in I : |(\on{S}_n)_{g(t)}| \text{ minimal}\} = \{t \in I : |\{g_1(t),\ldots,g_n(t)\}| \text{ maximal}\}
  \]
  which is open in $I$. 
  Choose $t_0 \in J$. There exists a permutation $\si_{t_0} \in \on{S}_n/(\on{S}_n)_{g(t_0)}$ so that $f(t_0)=\si_{t_0}.g(t_0)$. Set $\tilde g = \si_{t_0}.g$.
  We claim that $f=\tilde g$. This is true locally near $t_0$, since $J$ is open. 
  Assume for contradiction that there exists $t_1 \in I$ so that $f(t_1) \ne \tilde g(t_1)$. Without loss of generality assume $t_0 < t_1$ and let
  \[
  s = \sup\{t \in [t_0,t_1) : f|_{[t_0,t]}=\tilde g|_{[t_0,t]}\} \in (t_0,t_1).
  \] 
  But then the $\cC$-curve $h=f-\tilde g$ is identically $0$ on the (non-trivial) interval $[t_0,s]$ 
  and for each $\ep>0$ there exists $t \in (s,s+\ep)$ with $h(t) \ne 0$.
  Thus, $h$ must vanish of infinite order at $s$, which contradicts our assumption.  
\qed\end{demo}

\section{\texorpdfstring{$\cC$}{C}-polynomials} \label{secCpoly}

\subsection{Monic univariate complex polynomials}
The space of all monic univariate complex polynomials $P$ of fixed degree $n$, 
\begin{equation} \label{P}
P(z) = z^n + \sum_{j=1}^n (-1)^j a_j z^{n-j} = \prod_{j=1}^n (z-\la_j),\quad a_j,\la_j \in \C,
\end{equation}
naturally identifies with $\C^n$ (via $P \mapsto (a_1,\ldots,a_n)$).
It may also be
viewed as the orbit space $\C^n/\on{S}_n$ with respect to the standard action $\on{S}_n : \C^n$ of the
symmetric group $\on{S}_n$ on
$\C^n$ by permuting the coordinates (the roots $\la_j$ of $P$).
The elementary symmetric
functions 
\[
\si_j(\la_1,\ldots,\la_n) = \sum_{i_1<\cdots<i_j} \la_{i_1} \cdots \la_{i_j}
\] 
generate
the algebra of symmetric polynomials on $\C^n$, i.e., $\C[\C^n]^{\on{S}_n}=\C[\si_1,\ldots,\si_n]$.
It follows that the orbit projection $\C^n \to \C^n/\on{S}_n$ identifies with
the mapping $\si=(\si_1,\ldots,\si_n) : \C^n \to \C^n$ and we have $a_j = \si_j(\la_1,\ldots,\la_n)$ (Vieta's formulas).
The associated polynomials
\begin{equation} \label{De}
\De_k(\la_1,\ldots,\la_n) := \sum_{i_1 < i_2 < \cdots < i_k}
(\la_{i_1}-\la_{i_2})^2 \cdots (\la_{i_1}-\la_{i_k})^2 \cdots
(\la_{i_{k-1}}-\la_{i_k})^2 
\end{equation}
are symmetric. Thus there exist unique polynomials $\tilde \De_k$ such that $\De_k = \tilde \De_k \o \si$, 
and so the $\tilde \De_k$ are functions of $P$. 
The number of distinct roots of $P$ equals the maximal $k$ such that $\tilde \De_k(P) \ne 0$; it cannot decrease locally in $P$. 

If $P$ is any monic polynomial, we denote by $a_j(P)$ its coefficients so that $P$ takes the form \eqref{P} with $a_j=a_j(P)$. 

The inverse function property \thetag{\hyperref[C_5]{$\cC_5$}} and \thetag{\hyperref[C_1]{$\cC_1$}} imply the following lemma.

\begin{lemma}[Splitting lemma in $\cC$, see {\cite[3.2]{RainerQA}}] \label{split}
Let $P_0$ be a complex polynomial
satisfying
$P_0 = P_1 \cdot P_2$, where $P_1$ and $P_2$ are monic polynomials without common
root.
Then for $P$ near $P_0$ we have $P = P_1(P) \cdot P_2(P)$
for $\cC$-mappings
of monic polynomials $P \mapsto P_1(P)$ and $P \mapsto P_2(P)$,
defined for $P$
near $P_0$, with the given initial values.
(Here $P \mapsto P_i(P)$ is understood as a mapping $\R^{2n} \to \R^{2 \deg P_i}$.)
\end{lemma}

\subsection{\texorpdfstring{$\cC$}{C}-families of polynomials} \label{ssec:Cpoly}
By a \emph{$\cC$-family of polynomials} we mean a polynomial 
\begin{equation} \label{P(x)}
P(x)(z) = z^n + \sum_{j=1}^n (-1)^j a_j(x) z^{n-j},
\end{equation}
where the coefficients $a_j$ are complex-valued $\cC$-functions defined in a $\cC$-manifold $M$. 
Let $x_0 \in M$.
If $P(x_0)$ has distinct roots $\nu_1,\ldots,\nu_m$,  
the Splitting Lemma~\ref{split} provides a $\cC$-factorization
$P(x) = P_1(x) \cdots P_m(x)$ near $x_0$ such that no two factors have common roots and
all roots of $P_h(x_0)$ are
equal to $\nu_h$, for $1 \le h \le m$.
This factorization amounts
to a reduction of $\on{S}_n : \C^n$ to 
$\on{S}_{n_1} \times \cdots \times \on{S}_{n_m} :\C^{n_1} \oplus \cdots \oplus \C^{n_m}$, where $n_h$ is the
multiplicity of $\nu_h$. 
In this situation we shall write
\[
\on{S}(P(x_0)) := \on{S}_{n_1} \times \cdots \times \on{S}_{n_m}.
\]
In other words, $\on{S}(P(x_0))$ is the stabilizer of the ordered $n$-tuple consisting of the roots of $P(x_0)$ with multiplicities.  

Furthermore, we will \emph{remove fixed points} of 
$\on{S}_{n_1} \times \cdots \times \on{S}_{n_m} : \C^{n_1} \oplus \cdots \oplus \C^{n_m}$ or, equivalently, reduce each
factor $P_h$ to the
case $a_1(P_h)=0$ by replacing $z$
by $z-a_1(P_h)/n_h$. The effect on the roots of $P_h$ is a shift by a $\cC$-function.

For later reference we state the following result.

\begin{proposition}[{\cite[II Thm.\ 5.2]{Kato76}}] \label{controots}
  The roots of a polynomial \eqref{P(x)} with continuous coefficients $a_j : \R \to \C$ admit a continuous parameterization.
\end{proposition}

\subsection{Normal nonflatness} \label{norm-nf}

Let $I \subseteq \R$ be an open interval and let $I \ni t \mapsto P(t)$ be a $\cC$-family of polynomials \eqref{P(x)}.
We say that $P$ is \emph{normally nonflat at $t_0 \in I$} if it has the following property:
\begin{enumerate}
\item[\thetag{N}\phantomsection\label{N}] Let $k$ be maximal with the property that the germ at $t_0$ of $t \mapsto \tilde \De_k(P(t))$ is not $0$. 
Then $t \mapsto \tilde \De_k(P(t))$ is not infinitely flat at $t_0$.
\end{enumerate}
By \eqref{De}, condition \thetag{\hyperref[N]{N}} is equivalent to the following: 
Let $\la_j$ denote the germs at $t_0$ of a continuous parameterization of the roots of $P$; such exist by Proposition~\ref{controots}.
Then the order of contact at $t_0$ of any two unequal $\la_j$ is finite.
Evidently, \thetag{\hyperref[N]{N}} is satisfied if $P$ is a $\CQ$-polynomial.

We shall say that $P$ is \emph{normally nonflat} if \thetag{\hyperref[N]{N}} holds at each $t_0 \in I$.

\begin{lemma}[{\cite[2.1]{RainerAC}}] \label{mult}
  Let $P$ be a polynomial \eqref{P(x)} with coefficients $a_j : \R,0 \to \C$ germs at $0$ of $\cC$-functions, and $a_1=0$.
  Then, for integers $r$, the following conditions are equivalent:
  \begin{enumerate}
  \item $\om(a_j) \ge j r$, for all $2 \le j \le n$;
  \item $\om(\tilde{\Delta}_j) \ge j (j-1) r$, for all $2 \le j \le n$.
  \end{enumerate}
  Consequently,
  if $P$ is normally nonflat at $0$ and $\om(a_j)=\infty$ for all $j$, then $a_j=0$ for all $j$.
\end{lemma}

\begin{proposition}[Puiseux's theorem in $\cC$] \label{Puiseux}
  Let $P$ be a polynomial \eqref{P(x)} 
  with coefficients $a_j : \R,0 \to \C$ germs at $0$ of $\cC$-functions. 
  If $P$ is normally nonflat at $0$,
  then there exists a positive integer $\ga$ and germs $\la_j : \R,0 \to \C$ of $\cC$-functions such that
  $P(t^\ga)(z) = \prod_{j=1}^n (z- \la_j(t))$.
\end{proposition}

\begin{demo}{Proof}
  For $\cC=C^\infty$ this was proved in \cite[3.2]{RainerAC}. The same proof works for general $\cC$. See also \cite{LMRac}.  
\qed\end{demo}

\begin{lemma}[Glueing local choices of roots] \label{glue}
  Let $\R \ni t \mapsto P(t)$ be a $\cC$-curve of polynomials \eqref{P(x)}. 
  If $P$ is normally nonflat and locally admits $\cC$-parameterizations of its roots, i.e., for each 
  $t_0 \in \R$ there exist an open interval $I_{t_0} \ni t_0$ and $\cC$-functions which represent the roots of $P$ on $I_{t_0}$, 
  then there exists a global $\cC$-parameterization of the roots. 
\end{lemma}

\begin{demo}{Proof}
  Let $I \subseteq \R$ be a proper open subinterval and let $\la_j$, $1 \le j \le n$, be $\cC$-functions which represent 
  the roots of $P$ on $I$. We show that the $\cC$-parameterization $\la_j$ can be extended to a larger domain. 
  Let the right (say) endpoint $b$ of $I$ be finite. There exists a $\cC$-parameterization $\mu_j$ of the roots on some 
  open interval $I_b \ni b$. By Lemma~\ref{unique}, we may renumber the $\mu_j$ so that for all $j$, $\la_j = \mu_j$ on 
  their common domain $I \cap I_b$. So together the $\la_j$ and the $\mu_j$ form a $\cC$-parameterization of the roots on $I \cup I_b$.      
\qed\end{demo}

\section{\texorpdfstring{$\CQ$}{C_Q}-polynomials solvable along \texorpdfstring{$\CQ$}{C_Q}-arcs} \label{secpoly}

We have shown in \cite[6.7]{RainerQA} that a $\CQ$-polynomial $P$ admits a $\CQ$-parameterization of its roots after
desingularization by means of local blow-ups and local power substitutions. 
In this section we shall prove that local blow-ups suffice if $P$ is 
\emph{solvable along $\CQ$-arcs}. 
This will be applied to the characteristic polynomial of normal $\CQ$-matrices in Section~\ref{ramatrix}. 
It might also be of independent interest.

We say that a $\CQ$-family $M \ni x \mapsto P(x)$ of polynomials  
is \emph{solvable along $\CQ$-arcs} if, 
for all $\CQ$-curves $c : \R \to M$, 
the roots of $P \o c$ admit $\CQ$-parameterizations.

\begin{theorem} \label{poly}
Let $M$ be a $\CQ$-manifold and let $M \ni x \mapsto P(x)$ be a $\CQ$-family of polynomials \eqref{P(x)} solvable along $\CQ$-arcs.    
Let $K \subseteq M$ be compact.
Then there exists a finite covering $\{\pi_k : U_k \to W\}$ of a neighborhood $W$ of $K$, where each  
$\pi_k$ is a composite of finitely many local blow-ups, 
such that, for all $k$, the family of polynomials 
$P \o \pi_k$ allows a $\CQ$-parameterization of its roots on $U_k$.
\end{theorem}

\begin{demo}{Proof}
Since the statement is local, we may assume without loss of generality that $M$ is an open neighborhood of $0 \in \R^q$.
We use induction on the cardinality $|\on{S}(P(0))|$ of $\on{S}(P(0))$.

If $|\on{S}(P(0))|=1$, all roots of $P(0)$ are pairwise different.
So the statement follows from the $\CQ$-implicit function theorem \thetag{\hyperref[C_5]{$\cC_5$}} or from the Splitting Lemma~\ref{split}. 

Suppose that $|\on{S}(P(0))|>1$. Let $\nu_1,\ldots,\nu_m$ denote the distinct roots of $P(0)$; some of them are multiple ($m=1$ is allowed).
The Splitting Lemma~\ref{split} provides a $\CQ$-factorization
$P(x) = P_1(x) \cdots P_m(x)$ near $0$ such that the roots of distinct factors remain separated and $P_h(0)(z)= (z-\nu_h)^{n_h}$
for $1 \le h \le m$.
We reduce to $\on{S}_{n_1} \times \cdots \times \on{S}_{n_m} : \C^{n_1} \oplus \cdots \oplus \C^{n_m}$ 
and we remove fixed points (see \ref{ssec:Cpoly}), which preserves solvability along $\CQ$-arcs. 
So, if $a_{h,j} := a_j(P_h)$ denote the coefficients of $P_h$,
we may assume that $a_{h,1}=0$ for all $h$.
Then all roots of $P_h(0)$ are equal to $0$, and hence $a_{h,j}(0) = 0$, for all $1 \le h\le m$ and $1 \le j \le n_h$.
If all coefficients $a_{h,j}$ of $P_h$ are identically
$0$, so are all its roots, and we remove the factor $P_h$ from the product $P_1 \cdots P_m$.
Thus we can assume that for each $1 \le h \le m$ there is a $2 \le j \le n_h$ such that $a_{h,j} \ne 0$.

Let us define the $\CQ$-functions
\begin{equation} \label{Aa}
A_{h,j}(x) = a_{h,j}(x)^{\frac{n!}{j}} \quad (\text{for } 1 \le h \le m \text{ and } 2 \le j \le n_h).
\end{equation}
By Theorem~\ref{resth}, we find a finite covering $\{\pi_k : U_k \to U\}$ of a neighborhood $U$ of $0$ 
by $\CQ$-mappings $\pi_k$, each of
which is a composite of finitely many local blow-ups, such that, for each $k$, 
the non-zero $A_{h,j} \o \pi_k$ (for $1 \le h \le m$ and $2 \le j \le n_h$) and its pairwise non-zero differences 
$A_{h,i} \o \pi_k - A_{l,j} \o \pi_k$ (for $1 \le h \le l \le m$, $1 \le i \le n_h$, and $1 \le j \le n_l$) 
simultaneously have only normal crossings. 

Let $k$ be fixed and let $x_0 \in U_k$. 
Then $x_0$ admits a neighborhood $W_k$ with suitable coordinates in which $x_0=0$ and so that 
either $A_{h,j} \o \pi_k=0$ or 
\[
(A_{h,j} \o \pi_k)(x)=x^{\al_{h,j}} A_{h,j}^{k}(x), 
\]
where $A_{h,j}^{k}$ is a non-vanishing $\CQ$-function on $W_k$, and $\al_{h,j} \in \N^q$.
The collection of exponents $\{\al_{h,j} : A_{h,j} \o \pi_k \ne 0, 1 \le h \le m, 2 \le j \le n_h\}$ is totally ordered, by Lemma~\ref{order}. 
Let $\al$ denote its minimum. 

If $\al=0$, then $(A_{h,j} \o \pi_k)(x_0)=A_{h,j}^k(x_0)\ne 0$ for some $1 \le h \le m$ and $2 \le j \le n_h$. 
So, by \eqref{Aa}, not all roots of $(P_h \o \pi_k)(x_0)$ coincide (since $a_{h,1} \o \pi_k = 0$), and, thus, $|\on{S}((P \o \pi_k)(x_0))| < |\on{S}(P(0))|$.
Obviously, $P \o \pi_k$ is again solvable along $\CQ$-arcs.
By the induction hypothesis, 
there exists a finite covering $\{\pi_{kl} : W_{kl} \to W_k\}$ of $W_k$ (possibly shrinking $W_k$) of the required type such that, 
for all $l$, the family of polynomials $P \o \pi_k \o \pi_{kl}$ allows a $\CQ$-parameterization of its roots on $W_{kl}$. 

Let us assume that $\al \ne 0$.
Then there exist $\CQ$-functions $\tilde A_{h,j}^{k}$ on $W_k$ (maybe some of them $0$) such that, for all $1 \le h \le m$ and $2 \le j \le n_h$,
\begin{gather} \label{Atilde}
(A_{h,j} \o \pi_k)(x)=x^{\al} \tilde A_{h,j}^{k}(x), \text{ and},\\
\tilde A_{h,j}^{k} = A_{h,j}^{k} \text{ is non-vanishing, for some } 1 \le h \le m \text{ and } 2 \le j \le n_h. \label{Atilde2}
\end{gather}
Let us write
\[
\frac{\al}{n!} = \left(\frac{\al_1}{n!},\ldots,\frac{\al_q}{n!}\right) 
= \left(\frac{\be_1}{\ga_1},\ldots,\frac{\be_q}{\ga_q}\right), 
\]
where $\be_i,\ga_i \in \N$ are relatively prime (and $\ga_i>0$), for all $1 \le i \le q$. 

\begin{claim} \label{cl1}
 $\ga_i=1$ for all $1 \le i \le q$. 
\end{claim}

  We have to prove that $\al/n! \in \N^q$. Assume for contradiction that there is an $i_0$ such that $\al_{i_0}/n! \not\in \N$.
  Let $u \in W_k$ be such that $u_{i_0}=0$ and $u_i \ne 0$, for $i \ne i_0$, and let $e_{i_0}$ denote the $i_0$th standard unit vector in $\R^q$.
  Since $P_h$ is solvable along $\CQ$-arcs, we have 
  \[
  Q_h(t)(z):=P_h(\pi_k(u+te_{i_0}))(z)=\prod_{j=1}^{n_h} (z-\la_{h,j}(t))
  \] 
  for $\CQ$-functions $\la_{h,j}$ near $t=0$.
  By \eqref{Atilde2}, there exist $h_0$ and $2 \le j_0 \le n_{h_0}$ so that $\tilde A_{h_0,j_0}^{k}$ is non-vanishing.
  By \eqref{Aa} and \eqref{Atilde}, we have
  \begin{gather}
    \om(a_j(Q_{h_0})^{\frac{n!}{j}}) \ge \al_{i_0},  \text{ for all } j, \text{ and}, \label{a>}\\
    \om(a_{j_0}(Q_{h_0})^{\frac{n!}{j_0}}) = \al_{i_0}. \label{a=}
  \end{gather}
  Since $\al_{i_0}>0$, \eqref{a>} implies that
  $\la_{h_0,j}(0)=0$ for all $1 \le j \le n_{h_0}$. 
  Set 
  \[
  r_{h_0} := \min_{1 \le j \le n_{h_0}} \om(\la_{h_0,j}).
  \] 
  There exist $\CQ$-functions $\mu_{h_0,j}$ such that
  $\la_{h_0,j}(t) = t^{r_{h_0}}\mu_{h_0,j}(t)$ for all $j$, by \thetag{\hyperref[C_4]{$\cC_4$}}. 
  Then 
  \begin{equation} \label{r>}
  \om(a_j(Q_{h_0})) \ge j r_{h_0},  \text{ for all } j,
  \end{equation}
  and the $\mu_{h_0,j}$ parameterize the roots of the polynomial $\tilde Q_{h_0}$
  with coefficients $a_j(\tilde Q_{h_0}(t)) := t^{-j r_{h_0}}a_j(Q_{h_0}(t))$. Since $\mu_{h_0,j}(0) \ne 0$ for some $j$, 
  not all coefficients of $\tilde Q_{h_0}(0)$ vanish. 
  So, for some $j_1$, we have 
  \begin{equation} \label{r=}
    \om(a_{j_1}(Q_{h_0}))=j_1 r_{h_0}. 
  \end{equation} 
  Combining \eqref{a>} and \eqref{r=} we find $\al_{i_0}/n! \le r_0$, and \eqref{a=} and \eqref{r>} together imply  
  $\al_{i_0}/n! \ge r_0$.
  Hence $\al_{i_0}/n! = r_{h_0} \in \N$, a contradiction. Thus Claim~\ref{cl1} is shown.
\smallskip

By \eqref{Aa}, \eqref{Atilde}, and Claim~\ref{cl1}, each
$a_{h,j} \o \pi_k$ is divisible by $x^{j \be}$ where $\be=(\be_1,\ldots,\be_q)$, 
and, by \thetag{\hyperref[C_4]{$\cC_4$}}, 
there exist $\CQ$-functions $a_{h,j}^{k}$ on $W_k$ such that
\begin{equation} \label{P^k}
(a_{h,j} \o \pi_k)(x) = x^{j \be} a_{h,j}^{k}(x) \quad (\text{for } 1 \le h \le m \text{ and } 2 \le j \le n_h).
\end{equation}
Consider the $\CQ$-family of polynomials $P_h^k$ with coefficients $a_j(P_h^k) := a_{h,j}^{k}$.
By \eqref{Atilde2}, there exist $1 \le h \le m$ and $2 \le j \le n_h$ such that $a_{h,j}^{k}(x_0) \ne 0$,
and, hence, not all roots of $P_h^k(x_0)$ coincide.
So for $P^k := P_1^k \cdots P_m^k$ we have $|\on{S}(P^k(x_0))| < |\on{S}(P(0))|$.

\begin{claim} \label{cl2}
  $P^k$ is solvable along $\CQ$-arcs.
\end{claim}

Let $c : \R \to W_k$ be a $\CQ$-curve. 
By Lemma~\ref{glue}, it suffices to show that the roots of $P^k \o c$ locally admit $\cC$-parameterizations,
and without loss of generality it is enough to show this locally near $0 \in \R$.
By Proposition~\ref{Puiseux}, there exists $\ga \in \N_{>0}$ such that 
$t \mapsto P^k(c(t^\ga))$ admits a $\CQ$-parameterization $\la_j$ of its roots near $t=0 \in \R$.
Let $\ga$ be minimal with that property. 
For contradiction assume that $\ga >1$. 
By \eqref{P^k}, the roots of $P^k$ and $P \o \pi_k$ differ by the monomial factor $m(x):=x^{\be}$. 
Thus, the functions $\mu_j(t):= m(c(t^\ga)) \cdot \la_j(t)$ form a $\CQ$-parameterization of the roots of $t \mapsto P(\pi_k(c(t^\ga)))$.
Since $P \o \pi_k$ is solvable along $\CQ$-arcs, there exist $\CQ$-functions $\nu_j$ which parameterize the roots of $P \o \pi_k \o c$. 
Hence, both collections $\{\mu_j\}$ and $\{t \mapsto \nu_j(t^\ga)\}$ parametrize the roots of $t \mapsto P(\pi_k(c(t^\ga)))$, and, after renumbering, 
we may assume that $\nu_j(t^\ga)=m(c(t^\ga)) \la_j(t)$ for all $j$, by Lemma~\ref{unique}.
By \thetag{\hyperref[Q]{Q}} and \thetag{\hyperref[C_4]{$\cC_4$}}, the quotients $\nu_j/(m \o c)$ are $\CQ$-functions. 
As they parameterize the roots of $P^k \o c$, the choice of $\ga$ was not minimal, a contradiction.  
This proves Claim~\ref{cl2}.
\smallskip

Now, by the induction hypothesis, 
there exists a finite covering $\{\pi_{kl} : W_{kl} \to W_k\}$ of $W_k$ (possibly shrinking $W_k$) of the required type such that,
for all $l$, the family of polynomials $P^k \o \pi_{kl}$ admits a $\CQ$-parameterization 
$\la_{j}^{kl}$ of its roots on $W_{kl}$.
Then the $\CQ$-functions $x \mapsto m(\pi_{kl}(x)) \cdot \la_{j}^{kl}(x)$ form a choice of the roots 
of the family $x \mapsto (P \o \pi_k \o \pi_{kl})(x)$ for $x \in W_{kl}$.

Since $k$ and $x_0$ were arbitrary, the assertion of the theorem follows.
\qed\end{demo}

Let us call a $\CQ$-family $M \ni x \mapsto P(x)$ of polynomials \eqref{P(x)} \emph{solvable after blowing up} 
if the conclusion of Theorem~\ref{poly} holds, i.e.,
for $K \subseteq M$ compact, there exists a finite covering $\{\pi_k : U_k \to W\}$ of a neighborhood $W$ of $K$, where each  
$\pi_k$ is a composite of finitely many local blow-ups, 
such that, for all $k$, $P \o \pi_k$ allows a $\CQ$-parameterization of its roots.

\begin{corollary}[Solvability along \texorpdfstring{$\CQ$}{C_Q}-arcs and after blowing up are equivalent]
  A $\CQ$-family $M \ni x \mapsto P(x)$ of polynomials \eqref{P(x)} is solvable along $\CQ$-arcs if and only if 
  it is solvable after blowing up.
\end{corollary}

\begin{demo}{Proof}
 One direction is shown in Theorem~\ref{poly}. For the converse direction
 let $c : \R \to M$ be a $\CQ$-curve.
 By Lemma~\ref{glue}, it suffices to prove that $P \o c$ admits $\CQ$-parameterizations of its roots, locally.
 Let $t_0 \in \R$, 
 set $K := \{c(t_0)\}$, and apply the assumption that $P$ is solvable after blowing up.
 Denote by $c : \R,t_0 \to M$ the germ of $c$ at $t_0$.
 \[
 \xymatrix{
 && U_k \ar[rr] \ar[d]^{\pi_k} && \C^n \ar[d]^{\si} & \\
 \R,t_0 \ar[rr]_{c} \ar[rru] && W \ar[rr]_{P|_W} && \C^n/\on{S}_n \ar@{=}[r] & \C^n 
 }
 \]
 Since $\CQ$-curves admit a lifting over blow-ups,  
 all arrows in the diagram are of class $\CQ$.      
 This implies the statement.
\qed\end{demo}

\begin{remarks*}
  \thetag{1} Compare with the interrelation between \emph{arc-analyticity} and \emph{blow-analyticity}, see \cite{BM90} and \cite{Parusinski94}.
  
  \thetag{2} Hyperbolic $\CQ$-polynomials are solvable along $\CQ$-arcs, see \cite[6.11]{RainerQA}. 
  \emph{Hyperbolic} means all roots are real at each parameter value. 
  In the next section will meet another class of polynomials solvable along $\CQ$-arcs. 
\end{remarks*}

\section{Smooth perturbation theory for normal matrices} \label{ramatrix}

\begin{lemma} \label{fsolv}
  Let $P$ be a polynomial \eqref{P(x)} with coefficients $a_j : \R,0 \to \C$ germs at $0$ of $\cC$-functions,
  and assume that $P$ is normally nonflat at $0$.
  If there exist $\La_1,\ldots,\La_n \in \C[[t]]$ which represent the roots of the formal polynomial $\widehat P$, i.e., 
  \[
  \widehat P(t)(z) = z^n+\sum_{j=1}^n (-1)^j \widehat{a_j}(t) z^{n-j} = \prod_{j=1}^n (z-\La_j(t)),
  \]
  then there exist germs $\la_1,\ldots,\la_n : \R,0 \to \C$ of $\cC$-functions such that $P(t)(z) = \prod_{j=1}^n (z- \la_j(t))$ 
  and $\widehat{\la_j} = \La_j$ for all $j$.
\end{lemma}

\begin{demo}{Proof}
  In view of the reduction procedure described in \ref{ssec:Cpoly} (which preserves normal nonflatness) 
  we may assume that all roots of $P(0)$ equal $0$ and $a_1=0$.
  Let $r := \min_{1 \le j \le n} \om(\La_j) \ge 1$. 
  If $r=\infty$ then all $a_j = 0$, by Lemma~\ref{mult}, and by setting all $\la_j =0$ we are done.
  So we may assume that $r<\infty$. 
  For each $j$ we have $\om(\widehat{a_j})\ge j r$, 
  thus $a_j$ is divisible by $t^{j r}$, and, by \thetag{\hyperref[C_4]{$\cC_4$}}, there exist 
  $\cC$-germs $b_j$ such that $a_j(t)=t^{j r} b_j(t)$. 
  Consider the polynomial $Q$ with coefficients $a_j(Q):= b_j$. It is easy to see that $Q$ is normally nonflat at $0$ 
  and that not all roots of $Q(0)$ coincide. 
  Thus, induction on the cardinality of $\on{S}(P(0))$ proves the statement.
\qed\end{demo}

\begin{remark*}
  It is easy to check that the ring of germs at $0 \in \R$ of complex-valued $\CQ$-functions is a 
  Henselian excellent discrete valuation ring with maximal ideal $\mathfrak m = \{h : h(0)=0\}$ and 
  $\mathfrak m$-adic completion $\C[[t]]$. Thus, by \cite{Artin69}, \cite{Popescu86}, or \cite[Thm.\ 4.2]{Rotthaus87}, 
  it has the Artin approximation property which might be used alternatively to Lemma~\ref{fsolv} in the quasianalytic case. 
\end{remark*}

Let us introduce notation.
We associate with a parameterized family of complex matrices $A(x)=(A_{ij}(x))_{1 \le i,j \le n}$ its characteristic polynomial  
$\ch(A) := \det (A-z\I)$ and set $P_A := (-1)^n \ch(A)$.
Then $P_A$ is a family of polynomials \eqref{P(x)} with coefficients $a_j(P_A) = \on{Trace}(\La^j A)$, i.e.,
\begin{equation} \label{P_A}
P_A(x)(z) = (-1)^n \ch(A(x))(z) = z^n + \sum_{j=1}^n (-1)^j \on{Trace}(\La^j A(x)) z^{n-j}.
\end{equation}
We say that $A(x)=(A_{ij}(x))_{1 \le i,j \le n}$ is a family of \emph{normal} complex matrices if $A(x)A^*(x)=A^*(x)A(x)$ for all $x$.

\begin{proposition} \label{1dim}
  Let $A(t)=(A_{ij}(t))_{1 \le i,j \le n}$ be a $\cC$-curve of normal complex matrices, 
  i.e., the entries $A_{ij}$ belong to $\cC(\R,\C)$, such that $P_A$ is normally nonflat.   
  Then there exists a global $\cC$-parameterization of the eigenvalues and the eigenprojections 
  of $A$.
\end{proposition}

In the real analytic case the local statement of
this proposition is (by considering holomorphic extensions) a direct consequence of \cite[II Thm.\ 1.10]{Kato76} 
which exploits the monodromy of algebraic functions; see also \cite[3.5.1]{Baumgaertel85}.
An algebraic version for normal matrices over so-called \emph{Hermitian discrete valuation rings} is due to \cite{Adkins:1991ce}. 
Actually, for $\CQ$-curves of normal matrices, the local statement follows 
from \cite{Adkins:1991ce}, since the germs at $0 \in \R$ of complex-valued $\CQ$-functions 
form a Hermitian discrete valuation ring (as can be checked using Remark~\ref{fsolv}).

\begin{demo}{Proof}
  \emph{First we treat the eigenvalues.}
  By Lemma~\ref{glue},
  it suffices to show that there exist $\cC$-parameterizations of the eigenvalues, locally near each $t_0$. 
  Without loss of generality assume that $t_0=0$. In view of Lemma~\ref{fsolv} it is enough to prove the following claim.
  
  \begin{claim*}
    There exist $\La_1,\ldots,\La_n \in \C[[t]]$ such that 
    $\widehat{P_A}(t)(z) = \prod_{j=1}^n (z-\La_j(t))$.
  \end{claim*}
  
  This claim is a consequence of \cite{Adkins:1991ce}, since $\C[[t]]$ is a Hermitian discrete valuation ring and $\widehat{P_A} = P_{\widehat A}$,
  where the matrix $\widehat A(t)=(\widehat{A_{ij}}(t))$ is normal, since Taylor expansion commutes with transposition and conjugation
  (note that $\overline{\sum f_j t^j} = \sum \overline{f_j} t^j$). 
  
  Here is a direct proof more adapted to our situation.
  
  \begin{demo}{Proof of claim}
  Let $s$ be maximal with the property that the germ at $0$ of $\tilde \De_s(P_A)$ does not vanish identically.
  If $\tilde \De_s(P_A(0)) \ne 0$, then the Splitting Lemma~\ref{split} implies the assertion.
  So let us assume that $\tilde \De_s(P_A(0)) = 0$, i.e., generically distinct roots of $P_A$ meet at $0$.
  By Proposition~\ref{Puiseux}, there exists a minimal $\ga \in \N_{>0}$ such that  
  \begin{equation} \label{Pga}
  P_A(t^\ga)(z)=\prod_{j=1}^s(z-\la_j(t))^{m_j}
  \end{equation} 
  for generically distinct $\cC$-germs $\la_j : \R,0 \to \C$.
  Let $\th$ be a primitive $\ga$th root of unity and consider the formal power series 
  \begin{equation} \label{whla}
  \widehat{\la_j}(\th t) = \sum_{k\ge 0} \la_{j,k} \cdot (\th t)^k = \sum_{k\ge 0} \la_{j,k} \th^k \cdot t^k\quad (\text{where } \la_{j,k}=k!^{-1} \la_j^{(k)}(0)). 
  \end{equation}
  By \eqref{Pga},
  the $\widehat{\la_j}(\th t)$ represent the roots of the formal polynomial $\widehat{P_A}(t^\ga)$, likewise with the $\widehat{\la_j}(t)$. 
  Since $\C[[t,z]]$ is a unique factorization domain, we have:
  \begin{equation} \label{perm}
  \text{There exists $\si \in \on{S}_s$ such that $\widehat{\la_j}(\th t) = \widehat{\la_{\si(j)}}(t)$ for all $1 \le j \le s$.}  
  \end{equation} 
  We shall show that $\si$ is trivial. Then, in view of \eqref{whla}, 
  $\la_{j,k} \cdot \th^k = \la_{j,k}$ for all $j$ and all $k \in \N$. 
  So $\la_{j,k} = 0$ whenever $k \not\in \ga \N$, 
  and, thus, $\widehat{\la_j}(t^{1/\ga})$ is a formal power series in $t$. 
  By \eqref{Pga}, the formal power series $\widehat{\la_j}(t^{1/\ga})$, $1 \le j  \le s$, represent the distinct roots of $\widehat{P_A}(t)$;
  they are pairwise distinct by normal nonflatness. The claim follows.

  Suppose that $\si$ is non-trivial.
  Clearly, $\la_1,\ldots,\la_s$ parameterize the generically distinct eigenvalues of $t \mapsto A(t^\ga)$.
  Let $P_1,\ldots,P_s$ denote the respective eigenprojections:
  \begin{equation} \label{Sylvester}
  P_i(t) = \prod_{\substack{j=1 \\ j \ne i}}^s\frac{A(t^\ga)-\la_j(t)}{\la_i(t)-\la_j(t)}. 
  \end{equation}
  Normal nonflatness implies that 
  there exist (matrix-valued) $\cC$-germs $Q_i$ such that $P_i(t) = t^{-p_i} Q_i(t)$, $p_i \in \N$. 
  Since $A(t^\ga)$ is normal, and, thus, $\|P_i(t)\|=1$, each $P_i$ is of class $\cC$, by \thetag{\hyperref[C_4]{$\cC_4$}}.
  So we may consider the formal power series (with coefficients $n \times n$ matrices)
  $\widehat{P_i}(\th t) = \sum_{k \ge 0} P_{i,k} \cdot (\th t)^k = \sum_{k \ge 0} P_{i,k}  \th^k \cdot t^k$, and 
  \eqref{perm} and \eqref{Sylvester} imply 
  that
  $\widehat{P_i}(\th t) = \widehat{P_{\si(i)}}(t)$ for all $i$. 
  If $\si$ is non-trivial, we get in particular  $P_{i,0} = P_{j,0}$ for some $i \ne j$.
  The fact that $P_i(t)P_j(t) = 0$ off $0$ implies $P_{i,0} P_{j,0} = 0$ and,   
  since $P_i$ is idempotent, we have $(P_{i,0})^2=P_{i,0}$. Therefore,
  \[
  P_{i,0} = (P_{i,0})^2 = P_{i,0} P_{j,0} = 0,
  \]
  which contradicts $\|P_i(t)\| = 1$. Hence $\si=\on{id}$ and the claim is shown.
  \smallskip
  
  \emph{Now we treat the eigenprojections.}
  Let $\la_j : \R\to \C$, $1\le j\le s$, be a global $\cC$-parameterizations of the generically distinct 
  eigenvalues of $A$ and let $P_j$, $1\le j\le s$, be the respective eigenprojections.
  Then each $P_i$ is expressed by \eqref{Sylvester} with $\ga=1$, and 
  we may conclude similarly as above that each eigenprojection is globally of class $\cC$.
  Normal nonflatness implies that points where distinct eigenvalues meet cannot accumulate.
  \qed\end{demo}
\end{demo}

\begin{theorem} \label{matrix}
  Let $M$ be a $\CQ$-manifold and let $A(x)=(A_{ij}(x))_{1 \le i,j \le n}$ be a family of normal complex matrices with entries $A_{ij}$ in $\CQ(M,\C)$. 
  Let $K \subseteq M$ be compact.
  Then there exists a finite covering $\{\pi_k : U_k \to W\}$ of a neighborhood $W$ of $K$, where each  
  $\pi_k$ is a composite of finitely many local blow-ups,
  such that, for all $k$, the family of normal complex matrices
  $A \o \pi_k$ allows a $\CQ$-parameterization of its eigenvalues and its eigenvectors on $U_k$.
  
  If $M = \R$, $A$ is of class $\cC$, and $P_A$ is normally nonflat, 
  then there exist global $\cC$-parameterizations of the eigenvalues and local $\cC$-parameterizations of the eigenvectors of $A$.
  If we assume \thetag{\hyperref[C_6]{$\cC_6$}}, also the eigenvectors admit a global $\cC$-parameterization.  
\end{theorem}

\begin{demo}{Proof} 
The proof is subdivided into several claims.

\begin{claim} \label{cl3}
  The statements about the eigenvalues are true.
\end{claim}

For $M = \R$ this was shown in Proposition~\ref{1dim}. 
Let $M$ be a general $\CQ$-manifold. Proposition~\ref{1dim} implies that 
the associated $\CQ$-family of polynomials $P_A$ is solvable along $\CQ$-arcs.
So Theorem~\ref{poly} implies Claim~\ref{cl3}.

\begin{claim} \label{cl4}
  Let $A=A(x)$ be a family of normal complex $n \times n$ matrices, where the entries $A_{ij}$ are $\CQ$-functions 
  and the eigenvalues of $A$ admit a $\CQ$-parameterization $\la_j$ in a neighborhood of $0 \in \R^q$.
  Then there exists a finite covering $\{\pi_k : U_k \to U\}$ of a neighborhood $U$ of $0$, where each  
  $\pi_k$ is a composite of finitely many local blow-ups, such that, for all $k$, 
  $A \o \pi_k$ admits a $\CQ$-parameterization of its eigenvectors.
\end{claim}

We prove Claim~\ref{cl4} using induction on $|\on{S}(P_A(0))|$.

First consider the following reduction:
Let $\nu_1,\ldots,\nu_m$ denote
the pairwise distinct eigenvalues of $A(0)$
with respective multiplicities $n_1,\ldots,n_m$. 
The sets
\[
\La_h:=\{\la_i : \la_i(0)=\nu_h\}, \quad 1 \le h \le m, 
\]
form a partition of the $\la_i$ such that $\la_i(x) \ne \la_j(x)$, for $x$ near $0$, if $\la_i$ and $\la_j$ belong to different $\La_h$. 
Consider
\begin{align*}
V_x^{(h)} := \bigoplus_{\la \in \La_h}
\on{ker}(A(x)-\la(x)) 
= \on{ker} \big(\circ_{\la \in \La_h}
(A(x)-\la(x))\big), \quad 1 \le h \le m.
\end{align*}
(The order of the compositions is not
relevant.)
Then $V_x^{(h)}$ is the kernel of a $\CQ$-vector bundle homomorphism $B(x)$
with constant rank (even of constant dimension of the kernel),
and thus it is a $\CQ$-vector subbundle of the trivial bundle 
$U \times \C^n \to U$ (where $U\subseteq \R^q$ is a neighborhood of $0$) which admits a $\CQ$-framing.
This can be seen as follows: Choose a basis of $\C^n$ such that $A(0)$
is diagonal. By the elimination procedure
one can construct a basis for the kernel of $B(0)$. For $x$ near $0$, the
elimination procedure (with the same choices)
gives then a basis of the kernel of $B(x)$. 
This clearly involves only operations which preserve the class $\CQ$.
The elements of this basis are
then of class $\CQ$ in $x$ near $0$.

Therefore, it suffices to find $\CQ$-eigenvectors in each subbundle
$V^{(h)}$ separately, expanded in
the constructed $\CQ$-frame field. But in this frame field the vector
subbundle looks again like a constant vector space.
So we may treat each of these parts ($A$ restricted to $V^{(h)}$, as
matrix with respect to the frame field) separately. 
For simplicity of notation we suppress the index $h$.

Let us write $a_j := a_j(P_A)$.
Suppose that all eigenvalues of $A(0)$ coincide and are equal to $a_1(0)/n$, 
according to \eqref{P_A}.
Eigenvectors of $A(x)$ are also eigenvectors of $A(x) - (a_1(x)/n) \I$ (and vice versa),
thus we may replace $A(x)$ by $A(x) - (a_1(x)/n) \I$ and assume that $a_1=0$.  
So $A(0)=0$.

If $A=0$ identically, we choose the eigenvectors constant and we are done.
Note that this proves Claim~\ref{cl4}, if $|\on{S}(P_A(0))|=1$.

Assume that $A \ne 0$.
By Theorem~\ref{resth}, there exists a finite covering 
$\{\pi_k : U_k \to U\}$ 
of a neighborhood $U$ of $0$, where each $\pi_k$ is a 
composite of finitely many local blow-ups, 
such that, for each $k$,
the non-zero entries $A_{ij} \o \pi_k$ of $A \o \pi_k$ and its pairwise non-zero differences 
$A_{ij} \o \pi_k - A_{lm} \o \pi_k$ simultaneously have only normal crossings. 

Let $k$ be fixed and let $x_0 \in U_k$.
Then $x_0$ admits a neighborhood $W_k$ with suitable coordinates in which $x_0=0$ and such that 
either $A_{ij} \o \pi_k =0$ or 
\[
(A_{ij} \o \pi_k)(x) = x^{\al_{ij}} B_{ij}^k(x), 
\]
where $B_{ij}^k$ is a non-vanishing $\CQ$-function on $W_k$, and $\al_{ij} \in \N^q$.
The collection of exponents $\{\al_{ij} : A_{ij} \o \pi_k \ne 0\}$ 
is totally ordered, by Lemma~\ref{order}. Let $\al$ denote its minimum. 

If $\al=0$, then $(A_{ij} \o \pi_k)(x_0)=B_{ij}^k(x_0)\ne 0$ for some $1 \le i,j \le n$. 
Since $a_1 \o \pi_k = 0$, we may conclude that not all eigenvalues of 
$(A \o \pi_k)(x_0)$ coincide.
Thus, $|\on{S}(P_{A \o \pi_k}(x_0))| < |\on{S}(P_A(0))|$, and, by the induction hypothesis, 
there exists a finite covering $\{\pi_{kl} : W_{kl} \to W_k\}$ of $W_k$ (possibly shrinking $W_k$) of the required type such that, 
for all $l$, the family of normal matrices $A \o \pi_k \o \pi_{kl}$ allows a $\CQ$-parameterization of its eigenvectors on $W_{kl}$.

Assume that $\al \ne 0$.
Then there exist $\CQ$-functions $A_{ij}^k$ (maybe some of them $0$) such that, 
for all $1 \le i,j \le n$,
\[
(A_{ij} \o \pi_k)(x) = x^\al A_{ij}^k(x),
\]
and $A_{ij}^k(x) = B_{ij}^k(x) \ne 0$ for some $i,j$ and all $x \in W_k$.
So $A^k(x) = (A_{ij}^k(x))$ forms a $\CQ$-family of normal $n \times n$ matrices, and its eigenvalues differ from those of 
$(A \o \pi_k)(x)$ by a monomial factor $x^\al$ and admit a $\CQ$-parameterization. Indeed, 
the $\CQ$-functions $\la_j \o \pi_k$ parameterize the eigenvalues of $A \o \pi_k$ and are divisible by $x^\al$, 
otherwise $x \mapsto \la_j(\pi_k(x))/x^\al$ would be an unbounded root of a polynomial with bounded coefficients, 
a contradiction (see e.g.\ \cite[2.4]{RainerQA}). In view of \eqref{P_A}, the $\CQ$-functions $x \mapsto \la_j(\pi_k(x))/x^\al$ represent 
the eigenvalues of $A^k$.

Eigenvectors of $A^k(x)$ are also eigenvectors of $(A \o \pi_k)(x)$ (and vice versa). 
As $A_{ij}^k(x_0) \ne 0$ for some $i,j$ and since $a_1(P_{A^k}) = 0$, 
not all eigenvalues of $A^k(x_0)$ coincide. Hence, $|\on{S}(P_{A^k}(x_0)| < |\on{S}(P_A(0))|$, and 
the induction hypothesis implies the statement. The proof of Claim~\ref{cl4} is complete.

\begin{claim} \label{cl5}
  If $M = \R$, $A$ is of class $\cC$, and $P_A$ is normally nonflat,
  then there exist local $\cC$-parameterizations of the eigenvectors of $A$.
  If we assume \thetag{\hyperref[C_6]{$\cC_6$}}, there exists a global $\cC$-parameterization of the eigenvectors.
\end{claim}

By Claim~\ref{cl3} the eigenvalues admit global $\cC$-parameterizations $\la_j$ on $\R$, which are unique up to a 
constant permutation, by Lemma~\ref{unique}.
The proof of Claim~\ref{cl4} works in this case as well: Theorem~\ref{resth} and Lemma~\ref{order}, the only ingredients that 
need quasianalyticity, are both trivially true, and normal nonflatness is preserved by the reduction process. 
So there are local $\cC$-choices of the eigenvectors.
The proof of Claim~\ref{cl4} further gives us,
for each eigenvalue $\la_j: \R \to \C$ with generic multiplicity $n_j$, a unique $n_j$-dimensional 
$\cC$-vector subbundle $V^{(j)}_t$ of $\R \times \C^n$ whose fiber over $t \in \R$ consists of eigenvectors for the 
eigenvalue $\la_j(t)$. 
By Proposition~\ref{1dim}, the eigenprojection $P_j$ corresponding to $\la_j$ is $\cC$ on $\R$ and $P_j(t)(\C^n)=V^{(j)}_t$.
It suffices to prove that each $P_j$ has a \emph{transformation function} of class $\cC$, cf.\ \cite[II \S 4.2]{Kato76},
i.e., there exists a matrix-valued function $\R \ni t \mapsto U_j(t)$ such that $U_j^{-1}(t)$ is invertible for each $t$, 
both $U_j$ and $U_j^{-1}$ are $\cC$ on $\R$, and $U_j(t)P_j(0)U_j^{-1}(t)=P(t)$. If $\{v_i\}$ is a basis 
of $P_j(0)(\C^n)$, then $\{U_j(t)v_i\}$ is a basis of $P_j(t)(\C^n)$. 

We construct a transformation function of class $\cC$ following \cite[II \S 4.2]{Kato76}.
Let us suppress the index $j$. Differentiation of $P^2=P$ and applying this identity several times yields 
\[
P'=[Q,P] = QP-PQ, \quad\text{ where }\quad Q=[P',P]=P'P-PP'.
\] 
By \thetag{\hyperref[C_5]{$\cC_5$}}, $Q$ is of class $\cC$. 
By \thetag{\hyperref[C_6]{$\cC_6$}}, the linear ODE 
\begin{equation} \label{diffe}
X'=QX
\end{equation}
with initial condition $X(0)=\I$
has a unique global solution $X=U$. Similarly, 
\begin{equation} \label{diffe2}
Y'=-YQ
\end{equation} 
with initial condition $Y(0)=\I$
has a unique global solution $Y=V$. Now $(VU)'=V'U+VU'=-VQU+VQU=0$ implies that $VU$ is a constant
and, by the initial conditions, we find that $VU=\I$, thus $U^{-1}=V$.
Since $(PU)'=P'U+PU'=(P'+PQ)U=QPU$, $PU$ is a solution of \eqref{diffe} with initial condition $X(0)=P(0)$.
Since the general solution of \eqref{diffe} is $X(t)=U(t)X(0)$, we have $U(t)X(0)=P(t)U(t)$, hence 
$U(t)P(0)U^{-1}(t)=P(t)$. So $U$ is a transformation function for $P$.

Moreover, $U(t)$ is unitary for each $t$ and hence the eigenvectors may be chosen orthonormal.
This is seen as follows, cf.\ \cite[II \S 6.2]{Kato76}:
Normality of $A$ implies $P^*=P$ and $(P')^*=P'$, by differentiation. Thus $Q =[P',P]$ is skew-Hermitian, 
and, since $U$ solves \eqref{diffe}, we find
\[
(U^*)'=-U^*Q
\] 
i.e., $U^*$ solves \eqref{diffe2}. Uniqueness implies that $U^*=V=U^{-1}$.
\qed\end{demo}

\section{Lipschitz eigenvalues of normal matrices} \label{Lipmatrix}  

There is the following result.

\begin{theorem}[{\cite{BhatiaDavisMcIntosh83}, \cite[VII.4.1]{Bhatia97}}] \label{d-Lip}
Let $A,B$ be normal complex $n \times n$ matrices and
let $\la_j(A)$ and 
$\la_j(B)$, $1 \le j \le n$, denote the respective eigenvalues.
Then 
\[
\min_{\si \in \on{S}_n} \max_{1\le j \le n} |\la_j(A)-\la_{\si(j)}(B)| 
\le C \|A-B\|
\]
for a universal constant $C$ with $1 < C < 3$,
where $\|~\|$ is the operator norm.
\end{theorem}

In particular, the \emph{unordered} $n$-tuple of eigenvalues $\la(A)=(\la_1(A),\ldots,\la_n(A))$ 
is continuous (even Lipschitz) as a function of 
the normal matrix $A$. However, the single eigenvalues do in general not allow continuous parameterizations,
see Example~\ref{excont}. Continuous parameterizations exist if $A$ is Hermitian (e.g.\ ordering by size 
$\la_j(A) \le \la_{j+1}(A)$; see \cite[4.1]{AKLM98}) or if $A$ depends on a single real parameter (see Proposition~\ref{controots}).
We shall show in this section that, 
if $A$ depends on parameters locally in a Lipschitz way and admits continuous parameterizations $\la_j$ of its eigenvalues, 
then the $\la_j$ are locally Lipschitz. 
No such result is true for the eigenvectors, see Section~\ref{optimal}.

We will repeatedly use the following fact.

\begin{lemma}[{\cite[4.3]{KLMR05}}] \label{inter}
  Let $c : (a,b) \to X$ be a continuous curve in a compact metric space $X$. 
  The set of accumulation points of $c(t)$ as $t \to a^+$ is connected.
\end{lemma}

Let us start with the one parameter case.

\begin{proposition} \label{Lip2}
  Let $A(t)=(A_{ij}(t))_{1 \le i,j \le n}$ be a curve of normal complex matrices, where the entries $A_{ij} : \R \to \C$  
  are locally Lipschitz. Then the eigenvalues of $A$ admit a parameterization which is locally Lipschitz.
  Actually, any continuous parameterization of the eigenvalues of $A$ is locally Lipschitz.
\end{proposition}

\begin{demo}{Proof}
  Let $s \in \R$ be fixed.
  Let $z$ be an eigenvalue of $A(s)$ of multiplicity $m$. 
  We choose a simple closed $C^1$-curve $\ga$ in the resolvent set of 
  $A(s)$ enclosing only $z$ among all eigenvalues of $A(s)$. 
  By continuity, see Proposition~\ref{controots}, no eigenvalue 
  of $A(t)$ lies on $\ga$, for $t$ near $s$; see also Lemma~\ref{resolv} below. 
  Now,
  \begin{equation*}
  t\mapsto -\frac1{2\pi i}\int_\ga (A(t)-z)\i\;dz =: P(t,\ga) = P(t)
  \end{equation*}
  is a locally Lipschitz curve of projections onto the direct sum of all 
  eigenspaces corresponding to eigenvalues of $A(t)$ in the interior of $\ga$ with constant rank (cf.\ Section~\ref{unbounded}).
  For $t$ near $s$,
  there are equally many eigenvalues in the 
  interior of $\ga$, and, by Proposition~\ref{controots}, we may call them 
  $\la_j(t)$, for $1\le j\le m$, so that each $\la_j$ is continuous. 
  
  The image of $t \mapsto P(t,\ga)$ describes a locally Lipschitz vector subbundle of the trivial bundle $\R \times \C^n \to \R$. 
  For each $t$ choose an orthonormal system of eigenvectors $v_j(t)$ of $A(t)$ corresponding to the $\la_j(t)$. 
  They form a (not necessarily continuous) framing. By local triviality of the vector bundle, 
  for each $t$ near $s$ and each sequence $t_k \to t$ there is a subsequence (again denoted by $t_k$) such that 
  $v_j(t_k) \to w_j(t)$, where the $w_j(t)$ form an orthonormal system of eigenvectors of $A(t)|_{P(t)(\C^n)}$.
  Consider
  \begin{equation} \label{inn}
  \frac{A(t)-\la_j(t)}{t_k-t} v_j(t_k) + \frac{A(t_k)-A(t)}{t_k-t} v_j(t_k) - \frac{\la_j(t_k)-\la_j(t)}{t_k-t}v_j(t_k) = 0.
  \end{equation}
  
  Now assume that $A'(s)$ exists.
  For $t=s$ take the inner product of \eqref{inn} with each $w_i(s)$: 
  The first summand vanishes, since all $\la_j(s)$ coincide with $z$ and since  
  the $w_i(s)$ form also an orthonormal system of eigenvectors of $A(s)^*$ corresponding to the eigenvalue $\overline z$  
  (cf.\ \cite[I \S 6.9]{Kato76}).
  Letting $k \to \infty$, we find that the $w_i(s)$ are a basis of eigenvectors of $P(s)A'(s)|_{P(s)(\C^n)}$ with eigenvalues
  \[
  \lim_{k \to \infty} \frac{\la_i(t_k)-\la_i(s)}{t_k-s}.
  \] 
  We may conclude, by Lemma~\ref{inter}, that the right-sided derivative $\la_j^{(+)}(s)$ of each $\la_j$ exists at $s$.
  Similarly, the left-sided derivatives $\la_j^{(-)}(s)$ exist and they form the same set of numbers with correct multiplicities.  
  Hence, applying a suitable permutation on one side of $s$ provides a continuous choice of the eigenvalues through $z$ which is differentiable at $s$.

  If we take the inner product of \eqref{inn} with $w_j(t)$ (for $t$ near $s$) and proceed to the limit, 
  then (as the first summand vanishes again by the same reason) we obtain
  \begin{equation} \label{right}
     \la_j^{(+)}(t) = \< A'(t)w_j(t) \mid w_j(t) \> \quad \text{ whenever $A'(t)$ exists}, 
  \end{equation}
  for a unit eigenvector $w_j(t)$ of $A(t)$ with eigenvalue $\la_j(t)$.
  A similar formula holds for the left-sided derivatives $\la_j^{(-)}(t)$.
  
  An inspection of these arguments shows that they hold for any continuous parameterization $\la_j$ of the eigenvalues of $A$.
  Hence we have shown:
  
  \begin{claim} \label{cl6} 
    Let $\la_j$ be any continuous parameterization of the eigenvalues of $A$.
    If $A'(s)$ exists, then the one-sided derivatives of $\la_j$ exist at $s$, left- and right-sided derivatives form the same
    set with correct multiplicities, namely, the set of eigenvalues of $A'(s)$,  
    and they satisfy a formula of type \eqref{right}.
    Applying a suitable permutation on one side of $s$ provides a continuous choice of the eigenvalues which is differentiable at $s$.
  \end{claim}
   
  Next we claim that each $\la_j$ is locally absolutely continuous. Then $\la_j$ is differentiable almost everywhere and its derivative is 
  locally bounded, by \eqref{right}. Thus $\la_j$ is locally Lipschitz.
  
  \begin{claim} \label{cl7}
    Any continuous parameterization $\la_j$ of the eigenvalues of $A$ is locally absolutely continuous.
  \end{claim}
  
  Taking the inner product of \eqref{inn} with $w_j(t)$ leads to
  \begin{equation} \label{abs}
    \Big|\Big\<\frac{A(t_k)-A(t)}{t_k-t} v_j(t_k) \mid w_j(t)\Big\>\Big| 
    = \Big|\frac{\la_j(t_k)-\la_j(t)}{t_k-t}\Big| |\< v_j(t_k) \mid w_j(t)\>|,
  \end{equation}
  for unit eigenvectors $v_j(t_k)$, $w_j(t)$ of $A(t_k)$, $A(t)$ with eigenvalue $\la_j(t_k)$, $\la_j(t)$, respectively, 
  and such that $v_j(t_k) \to w_j(t)$.
    
  Let $I \subseteq \R$ be an open bounded interval, $J \supseteq \overline I$ an open neighborhood of the closure $\overline I$, 
  and let $C_J$ denote the Lipschitz constant of $A$ on $J$ 
  (with respect to the operator norm).
  If $t \in J$ and $J \ni t_k \to t$, $t_k \ne t$, then, after passing to a subsequence (again denoted by $t_k$) 
  so that $v_j(t_k) \to w_j(t)$, 
  there is, by \eqref{abs}, a $k_0=k_0(t,(t_k)) \in \N$ such that 
  \begin{equation} \label{lab}
  \Big|\frac{\la_j(t_k)-\la_j(t)}{t_k-t}\Big| \le 2 C_J,\quad \text{ for all } k \ge k_0. 
  \end{equation}
  
  Let $j$ be fixed. Consider the continuous functions 
  \[
  q_k(t) := \frac{\la_j(t+1/k)-\la_j(t)}{1/k}\quad \text{ and set }\quad C_k := \max_{t \in \overline I} |q_k(t)|.
  \]
  We claim that $C_k$ is bounded in $k$. Otherwise there exists a subsequence (again denoted by $C_k$) such that $C_k \nearrow \infty$.
  Choose $t_k \in \overline I$ such that $C_k=|q_k(t_k)|$. 
  Since $\overline I$ is compact, after passing to a subsequence, $t_k \to t_\infty \in \overline I$. 
  We may also assume that this convergence is \emph{fast}, i.e.,  
  for all $n \in \N$ the sequence $k^n(t_k-t_\infty)$ is bounded.
  If $t_k=t_\infty$ constantly, then $C_k = |q_k(t_\infty)| \le 2 C_J$ for sufficiently large $k$, by \eqref{lab}.
  So we may assume that $t_k \ne t_\infty$, and consider   
  \begin{align} \label{esCk}
    \begin{split}
      C_k &= |q_k(t_k)| = \Big|\frac{\la_j(t_k+1/k)-\la_j(t_k)}{1/k}\Big| \\
      &\le \Big|\frac{\la_j(t_k+1/k)-\la_j(t_\infty)}{t_k+1/k-t_\infty}\Big| \cdot (1 + k |t_k-t_\infty|) + 
      \Big|\frac{\la_j(t_k)-\la_j(t_\infty)}{t_k-t_\infty}\Big| \cdot k |t_k-t_\infty|.
    \end{split}
  \end{align}
  By \eqref{lab}, there is some $k_0 \in \N$ such that both difference quotients on the right hand side of \eqref{esCk} 
  are bounded by $2C_J$ for all $k \ge k_0$.
  (Here we pass first to a subsequence of $t_k$ and then in turn to a subsequence of $s_k:=t_k+1/k$, 
  and set $k_0 := \max\{k_0(t_\infty,(t_k)),k_0(t_\infty,(s_k))\}$.)
  This contradicts the assumption that $C_k$ is unbounded.

  Since $C_k = \max_{t\in \overline I} |q_k(t)|$ is bounded, the sequence of functions $q_k$ is bounded in $L^p(I)$, for any $p\ge 1$.
  Since $L^p(I)$ is reflexive if $1 < p < \infty$, for such $p$, there exists a subsequence (again denoted by the full sequence)
  and an element $\la_j' \in L^p(I)$ such that (see e.g.\ \cite[V Thm.\ 4.2]{Conway85})
  \[
  q_k = \Big(t \mapsto \frac{\la_j(t+1/k)-\la_j(t)}{1/k}\Big) \longrightarrow \la_j' \quad \text{ weakly in } L^p(I). 
  \]
  Thus, for a test function $\vh \in C_c^\infty(I)$,
  \begin{align*}
    \int_I \la_j' \vh dt &= \lim_{k \to \infty} \int_I \Big(\frac{\la_j(t+1/k)-\la_j(t)}{1/k}\Big) \vh(t) dt \\
    &= \lim_{k \to \infty} \int_I \la_j(t) \Big(\frac{\vh(t-1/k)-\vh(t)}{1/k}\Big) dt = - \int_I \la_j \vh' dt,
  \end{align*}
  where we used substitution and assumed that $k$ is sufficiently large so that $\on{supp}(\vh) \pm 1/k \subseteq I$.
  This shows that $\la_j'$ is the weak derivative of $\la_j$, and, hence, $\la_j \in W^{1,p}(I)$.
  It follows that there is an absolutely continuous function $\tilde \la_j$ on $I$ which coincides with $\la_j$ almost everywhere in $I$, 
  and, thus, on a dense subset of $I$.
  By continuity, $\la_j=\tilde \la_j$. The proof of Claim~\ref{cl7} is complete.
\qed\end{demo}

\begin{proposition} \label{2diff}
  Let $A(t)=(A_{ij}(t))_{1 \le i,j \le n}$ be a curve of normal complex matrices, where the entries $A_{ij} : \R \to \C$  
  are $C^1$ (resp.\ $C^2$). Then the eigenvalues of $A$ admit a parameterization which is $C^1$ (resp.\ twice differentiable).
\end{proposition}

\begin{demo}{Proof}
  The proof is subdivided into several claims. We use the notation in the proof of \ref{Lip2}.
  
  \begin{claim} \label{cl8}
    If $A$ is $C^1$, then the eigenvalues admit a $C^1$-parameterization.
  \end{claim}
  
  We use induction on $n$. Let $\la_j$ be a continuous parameterization of the eigenvalues of $A$ (see Proposition~\ref{controots}).
  If $s$ is such that not all $\la_j(s)$ coincide, then the set $\{1,\ldots,n\}$ decomposes into the
  subsets $\{j : \la_j(s)=w\}$, $w \in \C$. For $i$ and $j$ in different (non-empty) subsets, 
  we have $\la_i(t) \ne \la_j(t)$ for all $t$ in an open
  interval $I_s$ containing $s$. As in the proof of \ref{Lip2}, we may treat distinct subsets 
  separately (by considering $A(t)|_{P(t,\ga)(\C^n)}$, where $\ga$ encloses exactly one of the distinct eigenvalues of $A(s)$ at a time).
  By the induction hypothesis, Claim~\ref{cl8} holds on $I_s$.
  
  Let $I$ be an open interval containing only points $s$, where not all $\la_j(s)$ coincide. 
  Let $J\subseteq I$ be a maximal open subinterval on which Claim~\ref{cl8} holds.
  We claim that $J = I$. Otherwise an endpoint $a$ of $J$ belongs to $I$ and there is a $C^1$-parameterization of the eigenvalues
  on an open interval $I_a \ni a$. Choosing $s \in J \cap I_a$ and permuting one choice of eigenvalues on one side 
  of $s$ in a suitable way (see Claim~\ref{cl6}), we might extend the $C^1$-parameterization beyond $a$, 
  contradicting maximality of $J$.      
  
  The set $E$ of points, where all eigenvalues coincide, is closed, and on its complement (which is a disjoint union of open intervals)
  we may parameterize the eigenvalues by $C^1$-functions $\mu_j$.
  For each isolated point $s$ of $E$ we apply in turn the following arguments:
  Extending all $\mu_j$ to $s$ by the single $n$-fold eigenvalue of $A(s)$ provides a continuous parameterization near $s$.
  By Claim~\ref{cl6}, we may assume that the $\mu_j$ are differentiable at $s$ after applying a suitable permutation 
  to the right of $s$.
  We claim that the derivative of each $\mu_j$ is continuous at $s$. 
  Namely, let $t_k \to s$ and apply \eqref{right} to $t_k$,
  \begin{equation} \label{mu'}
    \mu_j'(t_k) = \< A'(t_k)w_j(t_k) \mid w_j(t_k) \>.
  \end{equation}
  Choose a subsequence such that
  the $w_j(t_k)$ converge. 
  Then \eqref{mu'} converges to one of the eigenvalues of $A'(s)$.
  We may conclude, by Lemma~\ref{inter}, that the limit $\lim_{t \to s^+} \mu_j'(t)$ exists and that it equals one of the eigenvalues 
  of $A'(s)$ (the same for $t \to s^-$).
  By the mean value theorem, for $\th \in (0,1)$,
  \[
  \mu_j'(s) = \lim_{h \to 0^\pm} \frac{\mu_j(s + h)-\mu_j(s)}{h} = \lim_{h \to 0^\pm} \mu_j'(s+\th h) = \lim_{t \to s^\pm} \mu_j'(t).
  \]
  
  Finally, we extend each $\mu_j$ by the single $n$-fold eigenvalues of $A(s)$ at each accumulation point $s$ of $E$.
  By Claim~\ref{cl6} and since $s$ is an accumulation point of $E$, all $\mu_j'(s)$ exist and coincide.
  Let $t_k \to s$. By \eqref{mu'}, the sequence $\mu_j'(t_k)$ is bounded, and, thus, has a convergent subsequence. 
  By passing to a subsequence again so that the $w_j(t_k)$ converge, we find, by \eqref{mu'}, that 
  $\mu_j'(t_k)$ converges to some eigenvalue of $A'(s)$.
  But the latter all coincide with $\mu_j'(s)$, by Claim~\ref{cl6}. 
  This implies that the $\mu_j'$ are continuous at $s$.
  The proof of Claim~\ref{cl8} is complete. 
  
  \begin{claim} \label{cl9}
    Assume that $A$ is $C^2$.
    For each $s$ there is a $C^1$-parameterization of the eigenvalues near $s$ which is twice differentiable at $s$.
  \end{claim}
  
  We may assume without loss of generality that $s=0$. By the usual reduction procedure 
  (i.e., treating distinct eigenvalues of $A(0)$ separately by restricting to $P(t,\ga)(\C^n)$ for suitable $\ga$ 
  and in turn replacing $A$ by $A-(a_1(P_A)/n) \I$)
  we may assume without loss of generality that $0$ is the only eigenvalue of $A(0)$.
  Then $A(t)=t \tilde A(t)$, where $t \mapsto \tilde A(t)$ is a $C^1$-curve of normal matrices. 
  By Claim~\ref{cl8}, there is a $C^1$-parameterization $\mu_j$ of the eigenvalues of $\tilde A$.
  Then the functions $t \mapsto t \mu_j(t)$ are twice differentiable at $0$ and represent the eigenvalues of $A$.

  \begin{claim} \label{cl10}
    If $A$ is $C^2$, then the eigenvalues of $A$ admit a parameterization which is twice differentiable at every point.
  \end{claim}
  
  We modify the proof of Claim~\ref{cl8} and just indicate the necessary changes.
  
  Let $I$ be an open interval containing only points $s$ so that not all eigenvalues of $A(s)$ coincide.
  We show that a twice differentiable parameterization, say $\mu_j$, of the eigenvalues on an open subinterval $J \subseteq I$ 
  can be extended to $I$. Let $a \in I$ denote the right, say, endpoint of $J$. By induction, there exists a  
  twice differentiable parameterization $\la_j$ of the eigenvalues on an open interval $I_a \ni a$.
  Choose $s \in J \cap I_a$ and let $t_k \to s$.
  For each $k$ there is a permutation $\si \in \on{S}_n$ such that $\mu_j(t_k) = \la_{\si(j)}(t_k)$ for all $j$. 
  By passing to subsequences in turn (and Claim~\ref{cl6}), 
  we can assume that $\si$ does not depend on $k$ and that also $\mu_j'(t_k) = \la_{\si(j)}'(t_k)$ 
  for all $j$.
  Then 
  \begin{equation} \label{eq2diff}
  \frac{\mu_j'(t_k) - \mu_j'(s)}{t_k-s} = \frac{\la_{\si(j)}'(t_k) - \la_{\si(j)}'(s)}{t_k-s},
  \end{equation}
  and, thus, $\mu_j''(s) = \la_{\si(j)}''(s)$ for all $j$.
  So we can extend $\mu_j$ beyond $a$.
   
  Let $E$ denote the set of points $s$ so that all eigenvalues of $A(s)$ coincide.
  The last paragraph implies the existence of a twice differentiable parameterization of the eigenvalues on the 
  complement of $E$.
  By the arguments in the proof of Claim~\ref{cl8}, we may 
  construct from it a $C^1$-parameterization $\mu_j$ on $\R$ which is twice differentiable on the complement of $E$.
  Let $s \in E$ and $t_k \to s$.
  Let $\la_j$ be the parameterization of the eigenvalues near $s$ provided by Claim~\ref{cl9}. 
  After passing to subsequences as above, we have \eqref{eq2diff}.
  
  Assume that $s$ is isolated in $E$. 
  As $\la_j$ is twice differentiable at $s$, we conclude, 
  by Lemma~\ref{inter}, that the left-sided and the right-sided second order derivatives of $\mu_j$ exist at $s$, and  
  they form the same set of numbers with correct multiplicities.
  By applying a permutation to the right of $s$, we obtain a twice differentiable parameterization of the eigenvalues near $s$.
  We treat all isolated points $s \in E$ in this way.
  
  If $s$ is an accumulation point of $E$, then all $\mu_j'(s)$ coincide. Let $E \ni t_k \to s$. 
  In view of \eqref{eq2diff} and by Lemma~\ref{inter}, 
  we find that the second order derivatives of the $\mu_j$ exist at $s$ and they all coincide,
  by considering second order difference quotients on points in $E$.
  The proof is complete.
\qed\end{demo}

The following is a modification of \cite[Lemma]{KM03} and can be shown in the same way. 
For convenience of the reader, we include a proof.

\begin{lemma}[Cf.\ {\cite{KM03}}] \label{add}
  Let $I$ be an interval, $n \le N$, and $\mu_1,\ldots,\mu_N,\la_1,\ldots,\la_n : I \to \C$ be continuous 
  (resp.\ $C^1$ or twice differentiable) such that
  $|\{j : \la_j(t) = z\}| \le |\{j : \mu_j(t) = z\}|$ for all $t \in I$ and $z \in \C$.
  Then there exist continuous (resp.\ $C^1$ or twice differentiable) functions $\la_{n+1},\ldots,\la_N : I \to \C$ such that 
  $|\{1 \le j \le N : \la_j(t) = z\}| = |\{j : \mu_j(t) = z\}|$ for all $t \in I$ and $z \in \C$.
\end{lemma}

\begin{demo}{Proof}
We use induction on $N$. Certainly, the assertion is true if $N=1$.

For $s \in I$ such that not all $\mu_j(s)$ coincide, 
the sets $\{\la_j\}$ and $\{\mu_j\}$ decompose into subsets so that elements of different subsets do not meet 
on an open interval $I_s$ containing $s$. By induction, the statement holds on $I_s$.
     
Suppose that for no $t \in I$ all $\mu_j(t)$ coincide.  
Let $J$ be a maximal open subinterval of $I$ for which the statement of the lemma is true with $\la_j^1$ for $j > n$.
We will show $J=I$.
If the right (say) endpoint $b$ of $J$ belongs to $I$, then the statement holds on an open interval $I_b \ni b$
with $\la_j^2$ for $j > n$. Choose $s \in J \cap I_b$. 
We claim that there is a permutation $\si$ so that
each $\la_j^1$ in $\{t \in J : t \le s\}$ can be extended by $\la_{\si(j)}^2$ in $\{t \in I_b : t \ge s\}$, 
contradicting maximality of $J$. 
Let $t_k \to s^-$. We have $\la_j^1(t_k) = \la_{\si(j)}^2(t_k)$ for a permutation $\si$ which depends on $k$.
By passing to a subsequence, we may assume that $\si$ is independent of $k$ which shows the claim in the continuous case.
For the $C^1$ and the twice differentiable case, we pass to a subsequence again in order to obtain 
$(\la_j^1)'(t_k) = (\la_{\si(j)}^2)'(t_k)$ and we use the arguments surrounding \eqref{eq2diff}.

Let $E$ denote the closed set of all points in $I$, where all $\mu_j$ coincide.       
The complement $I \setminus E$ is a disjoint union of open intervals, on each of which the lemma holds.
Extending the $\la_j$ to $s \in E$ by the unique value $\mu_j(s)$, provides a continuous extension to $I$.
For the $C^1$ and the twice differentiable case, we may renumber the $\la_j$ to the right of each isolated point $s \in E$ 
so that they fit together in a $C^1$ or twice differentiable way (by Lemma~\ref{inter}).
If $s$ is an accumulation point of $E$, then all derivatives $\mu_j'(s)=:\mu'(s)$ coincide. 
Thus, by Lemma~\ref{inter}, each $\la_j$ is differentiable at $s$ with $\la_j'(s)=\mu'(s)$, and 
$\la_j'$ is continuous at $s$:
\[
\la_j'(t) = \mu_{\si_t(j)}'(t) \to \mu'(s) =\la_j'(s), \quad \text{ as } t \to s.
\] 
If $\mu_j$ is twice differentiable at $s$, then all $\mu_j''(s)=:\mu''(s)$ coincide,
by considering second order difference quotients on points in $E$. By Lemma~\ref{inter}, we may conclude that each $\la_j$ is twice differentiable at $s$ 
with $\la_j''(s)=\mu''(s)$. 
\qed\end{demo}

\subsection{The class \texorpdfstring{$\cL^{\cC}$}{LC}} \label{classL}

Let $U$ be open in $\R^q$.
We denote by $\cL^{\cC}(U)$ the class of all complex-valued
functions $f$ with the following properties:
\begin{enumerate}
  \item[\thetag{$\cL_1$}\phantomsection\label{L_1}] $f$ is defined and of class $\cC$ on the complement $U \setminus E_{U,f}$ 
  of a closed set $E_{U,f}$ with $\cH^q(E_{U,f})=0$ and $\cH^{q-1}(E_{U,f})<\infty$.
  \item[\thetag{$\cL_2$}\phantomsection\label{L_2}] $f$ is bounded on $U \setminus E_{U,f}$.
  \item[\thetag{$\cL_3$}\phantomsection\label{L_3}] $\nabla f$ is bounded on $U \setminus E_{U,f}$.
\end{enumerate}
$\cH^q$ denotes the $q$-dimensional Hausdorff measure.

\begin{theorem} \label{Lipthm}
  Let $x \mapsto A(x)=(A_{ij}(x))_{1 \le i,j \le n}$
  be a parameterized family of normal complex 
  matrices. Then:
  \begin{enumerate}
    \item If $x \mapsto A(x)$ is $\CQ$ in $x \in U$, where $U$ is open in $\R^q$, 
    then for any compact $K \subseteq U$ there exists a relatively compact neighborhood $W$ of $K$ and a 
    parameterization $\la_i$ of the eigenvalues of $A$ on $W$ which belongs to $\cL^{\CQ}$, thus, also to $SBV$.
    More precisely, the classical gradient $\nabla \la_i(x)$ exists for all $x \in W \setminus E_{W,\la_i}$ and for those $x$ we have
    \[
    \|\nabla \la_i(x)\|_{\infty} = \max_j |\p_j \la_i(x)|  \le \sup_{y \in W} \|A'(y)\|< \infty,
    \] 
    where $\|~\|$ is the operator norm and $A'(x)=d A(x)$ the 
    Fr\'echet derivative.
    \item If $x \mapsto A(x)$ is $C^{0,1}$ in $x \in U$, where $U$ is $c^\infty$-open in a convenient vector space $E$, 
    then each continuous eigenvalue $\la : U \supseteq V \to \C$, $V$ $c^\infty$-open, of $A$ is $C^{0,1}$. 
    If $x_0 \in U \cap \overline V$ and $c : \R \to U$ is a $C^\infty$-curve with $c(0)=x_0$ and $c((0,1]) \subseteq V$, 
    then $\la \o c|_{(0,1]}$ is globally Lipschitz on $(0,1]$.
  \end{enumerate}     
\end{theorem}

\begin{demo}{Proof}
  \thetag{1}  
  By \cite[9.6]{RainerQA}, there exists a parameterization $\la_i$ of the eigenvalues of $A$ on $W$ which satisfies 
  \thetag{\hyperref[L_1]{$\cL_1$}} and \thetag{\hyperref[L_2]{$\cL_2$}} and such that $\nabla \la_i \in L^1(W)$; 
  in particular, each $\la_i$ belongs to $SBV$, also by \cite[9.6]{RainerQA}.
  For $x \in W \setminus E_{W,\la_i}$,  
  $t \in \R$ small, and $e_j$ the $j$th standard unit vector in $\R^q$, 
  the curve $t \mapsto \la_i(x+te_j)$ represents an eigenvalue of $t \mapsto A(x+te_j)$, and, by Claim~\ref{cl6}, we have 
  \[
  |\p_j \la_i(x)| = \Big|\frac{d}{d t}|_{t=0} \la_i(x+te_j)\Big| \le \|A'(x)\| 
  \] 
  which implies the statement.
    
  \thetag{2}
  Suppose that $\la : V \to \C$ is a continuous eigenvalue of $A$.
  Let $c : \R \to V$ be $C^\infty$. Then $\la \o c$ parameterizes an eigenvalue of the $C^{0,1}$-curve of normal matrices $A \o c$.
  By Lemma~\ref{add}, $\la \o c$ can be completed to a continuous parameterization of the eigenvalues of $A \o c$ which is locally 
  Lipschitz, by Proposition~\ref{Lip2}, and so $\la \o c$ is locally Lipschitz. 
  Since $c$ was arbitrary, we conclude that $\la : V \to \C$ is $C^{0,1}$ (see \ref{ssec:def}).
  
  Let $x_0 \in U \cap \overline V$ and let $c : \R \to U$ be $C^\infty$ with $c(0)=x_0$ and $c((0,1]) \subseteq V$. 
  We already know that $\la \o c|_{(0,1]}$ is locally Lipschitz. Its derivative exists a.e.\ and is bounded by the Lipschitz constant 
  of $A \o c|_{[0,1]}$ (with respect to the operator norm), by Claim~\ref{cl6}. The assertion follows.
\qed\end{demo}

\begin{remark}
  Suppose that $\R^q \supseteq U \ni x \mapsto A(x)=(A_{ij}(x))_{1 \le i,j \le n}$ is a $C^{0,1}$-family of normal complex 
  matrices. Then Claim~\ref{cl6} actually implies that, whenever a (one-sided) directional derivative of an eigenvalue of $A$ exists, 
  it is uniformly bounded on compact subsets of $U$. 
\end{remark}

\begin{question} \label{question}
  Let $\R^q \supseteq U \ni x \mapsto A(x)=(A_{ij}(x))_{1 \le i,j \le n}$
  be a $C^{0,1}$-family of normal complex matrices.
  Do the eigenvalues of $A$ admit a parameterization by $SBV_{\on{loc}}$-functions whose classical gradient exits a.e.\ and 
  is locally bounded? 
\end{question}

\begin{corollary}
  Let $H_n(\C)$ (resp.\ $SH_n(\C)$) denote the real vector space of $n \times n$ Hermitian (resp.\ skew-Hermitian) matrices.
  For each $A \in H_n(\C)$, let $\la^{\uparrow}(A) = (\la_1(A),\ldots,\la_n(A))$ be the $n$-tuple of increasingly ordered 
  eigenvalues of $A$, i.e., $\la_i(A) \le \la_{i+1}(A)$.
  For each $A \in SH_n(\C)$, set $\mu^{\uparrow}(A) := \la^{\uparrow}(-i A)$.
  Then:
  \begin{enumerate}
    \item Both mappings $\la^{\uparrow} : H_n(\C) \to \R^n$ and $\mu^{\uparrow} : SH_n(\C) \to i \R^n$ are globally Lipschitz; 
    more precisely, 
    \begin{align*}
    \max_{j} |\la^{\uparrow}_j(A)- \la^{\uparrow}_j(B)| &\le \|A-B\|\quad \text{ for all } A,B \in H_n(\C), \text{ and} \\
    \max_{j} |\mu^{\uparrow}_j(A)- \mu^{\uparrow}_j(B)| &\le \|A-B\|\quad \text{ for all } A,B \in SH_n(\C),
    \end{align*}
    where $\|~\|$ is the operator norm.
    \item Let $\al \in (0,1]$. If $x \mapsto A(x) \in H_n(\C)$ (resp.\ $\in SH_n(\C)$) is $C^{0,\al}$ in $x \in U$, 
    where $U$ is $c^\infty$-open in a convenient vector space $E$, 
          then $x \mapsto \la^{\uparrow}(A(x))$ (resp.\ $x \mapsto \mu^{\uparrow}(A(x))$) forms a 
          $C^{0,\al}$-parameterization of the eigenvalues. 
  \end{enumerate}     
\end{corollary}

Note that \thetag{1} is due to \cite{Weyl12}, see also \cite[III.2.6]{Bhatia97}.
Compare \thetag{2} with \cite{KMR}.

\begin{demo}{Proof}
\thetag{1} The mapping $\la^{\uparrow}$ is continuous, and, by Theorem~\ref{Lipthm}\thetag{2}, it is  
locally Lipschitz, thus, differentiable a.e. Let $A, B \in H_n(\C)$ with $\|B\|=1$. 
Then $\R \ni t \mapsto \la^{\uparrow}_j(A+tB)$, $1 \le j \le n$, forms a continuous 
parameterization of the eigenvalues of $t \mapsto A+tB$. If $\la^{\uparrow}_j$ is differentiable at $A$, then 
Claim~\ref{cl6} implies 
$|d_B \la^{\uparrow}_j(A)| \le 1$, and, thus, $\|d \la^{\uparrow}_j(A)\| \le 1$. 
It follows that $\la^{\uparrow}$ is globally Lipschitz with Lipschitz constant $1$.  
The statement for $\mu^{\uparrow}$ follows immediately from $\mu^{\uparrow}(A) = \la^{\uparrow}(-i A)$.

\thetag{2} follows from \thetag{1}.  
\qed\end{demo}

\section{Perturbation theory for unbounded normal operators} \label{unbounded}

Let $E \ni x \mapsto A(x)$ be a parameterized family of unbounded normal operators in a Hilbert space $H$
with common domain of definition $V$ and with compact resolvent.
The parameter domain $E$ is either $\mathbb R$, $\mathbb R^n$, or an infinite dimensional convenient vector space,
respectively (as specified in Theorem~\ref{main}). 
Let $M=(M_k)$ be log-convex and stable under derivations 
(see \thetag{\hyperref[M_1]{M$_1$}} and \thetag{\hyperref[M_2]{M$_2$}}). 
In the Beurling case $C^{(M)}$ we also assume $M_k^{1/k} \to \infty$, or, equivalently $C^{(M)} \supseteq C^\om$ 
(see \thetag{\hyperref[M_3]{M$_3$}}).
Let $\sC$ stand for $C^\infty$, $C^{[M]}$, or $C^{k,\al}$; remember that $C^{[M]}$ means $C^{(M)}$ or $C^{\{M\}}$.

\begin{lemma}[Resolvent lemma \cite{KMRp}] \label{resolv}
  If $x \mapsto A(x)$ is $\sC$ in $x$, then the resolvent 
  $(x,z)\mapsto (A(x)-z)\i\in L(H,H)$ is $\sC$ on its natural domain,
  the global resolvent set
  $
  \{(x,z)\in E\times \C: (A(x)-z):V\to H \text{ is invertible}\}
  $
  which is open (and even connected).
\end{lemma}

\begin{demo}{Proof}
  For $C^\infty$, $C^{\{M\}}$, with special $M=(M_k)$, and $C^{0,\al}$ this was proved in \cite{KMRp}.
  The same proof works for general $M=(M_k)$, $C^{[M]}$, and $C^{k,\al}$; for the latter even with the same references.
  So we just sketch the proof for $C^{[M]}$:
  By definition $x \mapsto \<A(x)u \mid v\>$ is $C^{[M]}$ for each $u \in V$ and $v \in H$ 
  and, thus, $x \mapsto A(x)u$ is of the same class as a mapping $E \to H$ for each $u \in V$ (see \ref{ssec:def}).
  
  The following claim was proved in \cite[Claim (1)]{KMRp} for $C^{0,\al}$.
  
\begin{claim} \label{cl14}
  For each $x$ consider the norm $\|u\|_x^2:=\|u\|^2+\|A(x)u\|^2$ on 
  $V$. Since $A(x)$ is closed, $(V,\|~\|_x)$ is a 
  Hilbert space with inner product 
  $\langle u \mid v\rangle_x:=\langle u \mid v\rangle+\langle A(x)u \mid A(x)v\rangle$. 
  All these norms $\|~\|_x$ on $V$ are equivalent, 
  locally uniformly in $x$.
  We then equip $V$ with one of the 
  equivalent Hilbert norms, say $\|~\|_0$, and have $A(x) \in L(V,H)$ 
  for all $x$.
\end{claim}  

  By the linear uniform boundedness theorem and by \cite[5.1]{KMRu}, 
  we conclude that the mapping $E \to L(V,H), x \mapsto A(x)$, 
  is $C^{[M]}$. 
  If for some $(x,z) \in E \times \C$ the bounded operator $A(x)-z : V \to H$ is invertible, then this is true 
  locally with respect to the $c^\infty$-topology on the product which is the product topology, by \cite[4.16]{KM97}. 
  The resolvent $(x,z) \mapsto (A(x)-z)^{-1} : H \to V$ is $C^{[M]}$, since inversion is real analytic on the Banach 
  space $L(V,H)$ and since $C^{[M]} \supseteq C^\om$ is stable under composition \cite[4.11]{KMRu}.
\qed\end{demo}

\begin{demo}{Proof of Theorem~\ref{main}}
  Let $x_0 \in E$ and let $z$ be an eigenvalue of $A(x_0)$ of multiplicity $N$. 
  We choose a simple closed $C^1$-curve $\ga$ in the resolvent set of 
  $A(x_0)$ enclosing only $z$ among all eigenvalues of $A(x_0)$. 
  Since the global resolvent set is open, see Lemma~\ref{resolv}, no eigenvalue 
  of $A(x)$ lies on $\ga$, for $x$ near $x_0$.
  By Lemma~\ref{resolv}, 
  \begin{equation*}
  x\mapsto -\frac1{2\pi i}\int_\ga (A(x)-z)\i\;dz =: P(x,\ga) = P(x)
  \end{equation*}
  is a $\sC$-mapping. Each $P(x)$ is a projection, namely onto the direct sum of all 
  eigenspaces corresponding to eigenvalues of $A(x)$ in the interior of $\ga$,  with finite rank.
  Thus the rank must be constant: It is easy to see 
  that the  
  (finite) rank cannot fall locally, and it cannot increase, since the 
  distance in $L(H,H)$ of $P(x)$ to the subset of operators of 
  rank $\le N=\operatorname{rank}(P(x_0))$ is continuous in $x$ and is either 
  $0$ or $1$.   
  So, for $x$ in a neighborhood $U$ of $x_0$,
  there are equally many eigenvalues in the 
  interior of $\ga$, and we may call them 
  $\la_j(x)$ for $1\le j\le N$ (repeated 
  with multiplicity).

  The family of $N$-dimensional complex vector spaces 
  $U \ni x\mapsto P(x)(H)\subseteq H$ forms a $\sC$ Hermitian vector subbundle over $U$ 
  of the trivial bundle $U\times H\to U$:
  For given $x$, choose $v_1,\dots v_N\in H$ such that the $P(x)(v_i)$ are linearly independent and thus span $P(x)(H)$. 
  This remains true locally in $x$. We use the Gram Schmidt orthonormalization 
  procedure (which is $C^\om$ and preserves $\sC$) 
  for the $P(x)(v_i)$ to obtain a local orthonormal $\sC$-frame of the bundle. 

  Now $A(x)$ maps $P(x)(H)$ to itself and in a local $\sC$-frame it is given by a normal $N\times N$ matrix 
  parameterized in a $\sC$-way by $x$. 
  Then all local assertions (i.e., in a product neighborhood of $(x_0,z)$) of the theorem follow:
  \thetag{\hyperref[A]{A}} and \thetag{\hyperref[B]{B}} follow from Theorem~\ref{matrix}, 
  \thetag{\hyperref[C]{C}} and \thetag{\hyperref[D]{D}} from Theorem~\ref{Lipthm}, 
  \thetag{\hyperref[E]{E}} and \thetag{\hyperref[F]{F}} from Proposition~\ref{2diff}.
  
  Let us prove \thetag{\hyperref[D]{D}}. Let $E \subseteq U \ni x \to \la(x)$ be a continuous eigenvalue of
  $x \mapsto A(x)$ which is $C^{0,1}$ in $x \in E$, where $U$ is $c^\infty$-open in a convenient vector space $E$, 
  and let $c : \R \to U$ be a $C^\infty$-curve. We first show that $\la \o c$ is locally Lipschitz. 
  Let $t \in \R$ and $x=c(t) \in U$.
  By the local result, $x \in U$ has an open neighborhood $V$ such that 
  the restriction $\la|_{V}$ is $C^{0,1}$. Thus $\la|_V \o c|_I$ is locally Lipschitz, where $I$ is the connected 
  component of $c^{-1}(V)$ which contains $t$. This implies the statement. 
  
  For the supplements in \thetag{\hyperref[D]{D}} we need the following claim.
  
  \begin{claim} \label{cl15}
    Let $t\to A(t)$ be $C^{0,1}$ in $t\in \R$, let $I \subseteq \R$ be a compact interval, and
    let $t \mapsto \la_j(t)$ be a Lipschitz eigenvalue of $t \mapsto A(t)$ defined on a subinterval of $I$. 
    Then 
    \[
    |\la_j(s) - \la_j(t)| \le (1+|\la_j(t)|) (e^{C |s-t|}-1),
    \] 
    for a constant $C$ depending only on $I$.
  \end{claim}

  By reducing to $P(t)A(t)|_{P(t)(H)}$ as above, we may conclude that \eqref{abs} holds true, and, thus, 
  for $V_t=(V,\|~\|_t)$ and $\|u\|_t^2=\|u\|_H^2+\|A(t)u\|_H^2$,
  \begin{align*}
      \Big|\frac{\la_j(t_k)-\la_j(t)}{t_k-t}\Big| &|\< v_j(t_k) \mid w_j(t)\>|
      \le 
      \Big\|\frac{A(t_k)-A(t)}{t_k-t}\Big\|_{L(V_{t_k},H)} \|v_j(t_k)\|_{V_{t_k}} \|w_j(t)\|_H\\
      &= \Big\|\frac{A(t_k)-A(t)}{t_k-t}\Big\|_{L(V_{t_k},H)} \sqrt{\|v_j(t_k)\|_H^2 + \|A(t_k)v_j(t_k)\|_H^2}  \cdot 1\\
      &= \Big\|\frac{A(t_k)-A(t)}{t_k-t}\Big\|_{L(V_{t_k},H)} \sqrt{1 + |\la_j(t_k)|^2} \\
      &\le C (1 + |\la_j(t_k)|),  
  \end{align*}
  for a constant $C$, since all norms $\|~\|_t$ are uniformly equivalent locally in $t$, by Claim~\ref{cl14}.
  Since $t \mapsto \la_j(t)$ is Lipschitz, in particular, absolutely continuous, we obtain
  \[
  |\la_j'(t)| \le C+C |\la_j(t)| \quad \text{ a.e.,}
  \]
  and 
  Gronwall's lemma (e.g.\ \cite[(10.5.1.3)]{Dieudonne60}) implies the asserted inequality.
  \medskip  
  
  For the first supplement in \thetag{\hyperref[D]{D}},
  let $x_0 \in E \cap \overline U$, $c : \R \to E$, $c(0)=x_0$, and $c((0,1]) \subseteq U$.
  The continuous function $\la \o c|_{(0,1]}$ represents an eigenvalue of $A \o c|_{(0,1]}$ and is locally Lipschitz.
  By Claim~\ref{cl15}, $\la \o c|_{(0,1]}$ is bounded on $(0,1]$, and, by Lemma~\ref{inter}, 
  the limit $\lim_{t\to 0^+} \la \o c|_{(0,1]}(t)=:z$ exists and is an eigenvalue of $A(x_0)$. 
  The local result (the supplement in Theorem~\ref{Lipthm}\thetag{2}) yields that $\la \o c|_{(0,1]}$ is globally Lipschitz. 
  
  Finally, it remains to extend the local choices to global ones for the cases \thetag{\hyperref[A]{A}}, \thetag{\hyperref[D]{D}} if $E=\R$, 
  \thetag{\hyperref[E]{E}}, and \thetag{\hyperref[F]{F}}:

  \begin{claim} \label{cl13}
    In case \thetag{\hyperref[A]{A}} the eigenvalues and the eigenvectors admit global $C^\infty$ (resp.\ $C^{[M]}$) parameterizations.
    In the cases \thetag{\hyperref[D]{D}} if $E=\R$, \thetag{\hyperref[E]{E}}, and \thetag{\hyperref[F]{F}} the eigenvalues admit 
    global $C^{0,1}$, $C^1$, and   
    twice differentiable parameterizations, respectively.
  \end{claim}
  
  \emph{First we treat the eigenvalues.} 
  Let $\sE$ stand for $C^\infty$, $C^{[M]}$, $C^{0,1}$, $C^{1}$, or ``twice differentiable''; according to case    
  \thetag{\hyperref[A]{A}}, \thetag{\hyperref[D]{D}}, \thetag{\hyperref[E]{E}}, or \thetag{\hyperref[F]{F}}.
  Choose a numbering of the eigenvalues of $A(0)$ (with multiplicities).
  
  We consider sequences of $\sE$-functions $(\la_j)_{j \in \al}$, indexed by ordinals $\al$ 
  and defined on open intervals $I_j$ containing some fixed $t_0 \in \R$, which parameterize eigenvalues of $A$. 
  The set of all such sequences is partially ordered by inclusion of ordinals and then by restriction of the component functions.
  For each increasing chain the union is again such a sequence.
  By Zorn's lemma there exists a maximal sequence $(\la_j)$.
  
  In any maximal sequence each component function $\la_j$ is globally defined on $\R$. This is seen as follows:
  If $b<\infty$ is the right (say) endpoint of $I_j$, then, by Claim~\ref{cl15} and by Lemma~\ref{inter}, 
  the limit $\lim_{t \to b^-} \la_j(t) =: z$ exists and is an eigenvalue of $A(b)$.
  By the local results, there exist $\de,\ep>0$ such that all eigenvalues $|\la-z| < \ep$ of $A(t)$ for $|t-b|< \de$ 
  admit a parameterization by $\sE$-functions 
  \[
  \mu_j : (b-\de,b+\de) \to \{\la \in \C : |\la-z|<\ep\}=:B_\ep(z).
  \]
  In case \thetag{\hyperref[A]{A}}, $\la_j$ coincides with some $\mu_j$ on their common domain, 
  since unequal eigenvalues have finite order of contact,
  and, hence, it admits an extension beyond $b$.
  In the other cases 
  consider the $\la_j$ whose graph $\{(t,\la_j(t) : t \in I_j\}$ has non-empty intersection with the vertical boundary 
  $\{b-\de,b+\de\} \times B_\ep(z)$ 
  of the tube $(b-\de,b+\de) \times B_\ep(z) \subseteq \R \times \C$. 
  The endpoints of the corresponding intervals $I_j$ decompose $(b-\de,b+\de)$ into finitely many subintervals.
  We apply Lemma~\ref{add} on each subinterval; in case \thetag{\hyperref[D]{D}} where $\la_j \in C^{0,1}$ we use its 
  continuous version. Then we glue at the endpoints of the subintervals in a continuous, $C^1$, or twice differentiable way, respectively, 
  (as before in the proof of Proposition~\ref{2diff}) to obtain an extension 
  of at least $\la_j$. In case \thetag{\hyperref[D]{D}} this extension is $C^{0,1}$, since we already know that any 
  continuous eigenvalue is $C^{0,1}$.
  So the sequence was not maximal and the assertion follows.
  
  Any maximal sequence $(\la_j)$ parameterizes all eigenvalues of $A$ with the right multiplicities.
  If not, there is some $t_0$ and some eigenvalue $z$ of $A(t_0)$ such that $|\{j : \la_j(t_0) = z\}|$ is less than the multiplicity of $z$.
  By the local results, Lemma~\ref{add}, and the assumption on the order of contact in case \thetag{\hyperref[A]{A}}, 
  we may again conclude that $(\la_j)$ was not maximal, a contradiction. 
  
  \emph{Now let us treat the eigenvectors.}
  Let $\la_j : \R \to \C$ be a $C^\infty$ (resp.\ $C^{[M]}$) eigenvalue with generic multiplicity $N$. 
  By the arguments in the proof of Claim~\ref{cl5}, we obtain a unique global $N$-dimensional 
  $C^\infty$ (resp.\ $C^{[M]}$) vector subbundle of $\mathbb R\times H \to \R$ whose fiber over $t$ consists of eigenvectors for the 
  eigenvalue $\la_j(t)$. The corresponding $C^\infty$ (resp.\ $C^{[M]}$) eigenprojection $P_j$ 
  has a transformation function, since the arguments at the end of \ref{cl5} work in Banach spaces, see \cite{Yamanaka91} 
  and \cite[3.4]{Schindl09}. 
  So we find global $C^\infty$ (resp.\ $C^{[M]}$) eigenvectors for each eigenvalue. 
  This completes the proof of Claim~\ref{cl13} and the proof of the theorem.  
\qed\end{demo}

\begin{remark}[m-sectorial operators] \label{m-sec}
  The assumptions in Theorem~\ref{main} may be slightly relaxed, 
  if all $A(x)$ are m-sectorial operators. 
  In that case it suffices to assume that the associated quadratic forms $\fa(x)$ 
  have common domain of definition $V$ and $x \mapsto \fa(x)(u)$ is of the respective class for each $u \in V$. 
  In the following discussion we use the 
  definitions of \cite[VI]{Kato76}.  
     
  Let $E \ni x \mapsto \fa(x)$ be a parameterized family of closed sectorial (possibly unbounded) sesquilinear forms
  in a Hilbert space $H$ so that there is a dense subspace $V$ of $H$ which is the domain of definition of each $\fa(x)$, 
  i.e., $V(\fa(x))=V$. 
  We say that $\fa(x)$ is $C^\infty$, $C^{[M]}$, or $C^{k,\al}$ if $x \mapsto \fa(x)(u,v)$ is $C^\infty$, $C^{[M]}$, or $C^{k,\al}$
  for each $u,v \in V$; by polarization it is actually enough to require that $x \mapsto \fa(x)(u) = \fa(x)(u,u)$ 
  is of the respective class for all $u \in V$. Let $\sC$ stand for $C^\infty$, $C^{[M]}$, or $C^{k,\al}$.
  
  There is a bijective correspondence $\fa \mapsto A_\fa$ between the set of all densely defined closed sectorial forms $\fa$ and 
  the set of all m-sectorial operators $A$, where $\fa$ is bounded if and only if $A_\fa$ is bounded and $\fa$ is symmetric 
  (i.e., $\fa(u,v) = \overline{\fa(v,u)}$ for $u,v \in V$)
  if and only if $A_\fa$ is self-adjoint (by the first representation theorem \cite[VI Thm.~2.1]{Kato76}). 
  Note that an m-sectorial operator necessarily is densely defined and closed.
  
  Thus with the $\sC$-family $x \mapsto \fa(x)$ of closed sectorial forms we associate the family $x \mapsto A_\fa(x)=A_{\fa(x)}$ 
  of m-sectorial operators. If we also assume that $A_\fa(x)$ is normal for every $x$ and has compact resolvent for every 
  (equivalently, some) $x$, 
  then the conclusions of Theorem~\ref{main} hold true for the family $x \mapsto A_\fa(x)$. This follows from the 
  following two claims which replace Lemma~\ref{resolv} and Claim~\ref{cl15}.    
  
  \begin{claim} \label{clnew1}
    The mapping $(x,z) \mapsto (A_\fa(x)-z)^{-1} \in L(H,H)$ is $\sC$.
  \end{claim}

  This claim can be shown along the lines of the proof of \cite[VII Thm.\ 4.2]{Kato76}:
  Fix $x_0$.
  Without loss of generality $\fs = \on{Re} \fa(x_0) \ge 1$; this can be achieved by adding a suitable constant to $\fa(x_0)$. 
  Then the associated operator $S = A_{\fs} \ge 1$ is self-adjoint 
  and has a unique square root $G = S^{1/2}$. 
  Consider the forms 
  \[
  \fb(x)(u,v) = \fa(x)(G^{-1}u,G^{-1}v).
  \] 
  Each form $\fb(x)$ is defined everywhere on $H$, since $G^{-1}u \in V(G)=V(\fs)=V$ 
  (by the second representation theorem \cite[VI Thm.\ 2.23]{Kato76}), closable and thus bounded.
  The assumption that $x \mapsto \fa(x)$ is a $\sC$-family immediately gives that $x \mapsto \fb(x)(u,v)$ is $\sC$ for each 
  $u,v \in H$. 
  Consider the family of operators $B(x) \in L(H,H)$ defined by
  \[
  \< B(x)u \mid v \> = \fb(x)(u,v).
  \]    
  By the linear uniform boundedness principle and the fact that it suffices to use a set of linear functionals which together recognize 
  bounded sets instead of the whole dual (see the references in \ref{ssec:def}), $x \mapsto B(x) \in L(H,H)$ is $\sC$. 
  Replacing $u,v$ by $Gu,Gv$ we obtain
  \begin{equation} \label{faB}
  \fa(x)(u,v) = \< B(x)Gu \mid Gv \>, \quad \text{ for } u,v \in V. 
  \end{equation}
  So we have
  \[
  \< A_\fa(x)u \mid v \> = \< B(x)Gu \mid Gv \>, \quad \text{ for } u \in V(A_\fa(x)), v \in V,
  \]
  whence $G B(x) Gu$ exists and equals $A_\fa(x)u$, since $G$ is self-adjoint. 
  Since $G B(x) G$ is accretive and $A_\fa(x)$ is m-accretive, we have 
  \begin{align*}
  A_\fa(x) &= G B(x) G, \quad \text{ and } \\
  A_\fa(x)^{-1} &= G^{-1} B(x)^{-1} G^{-1} \quad \text{ near } x_0,
  \end{align*}
  where $G^{-1} \in L(H,H)$ and $B(x)^{-1} \in L(H,H)$, since $B(x_0)$ is invertible (cf.\ \cite[VI Thm.\ 3.2]{Kato76}).
  It follows that $x \mapsto A_\fa(x)^{-1} \in L(H,H)$ is $\sC$ near $x_0$. 
  Here we use that $\sC$ is preserved by composition with a real analytic mapping.
  
  Assume that $A_\fa(x_0)-z_0$ is invertible. Then $A_\fa(x)-z$ is invertible for $(x,z)$ near $(x_0,z_0)$.   
  For such $(x,z)$ we have
  \[
  (A_\fa(x)-z) A_\fa(x)^{-1} = 1 - z A_\fa(x)^{-1}
  \]
  and $1 - z A_\fa(x)^{-1} : H \to H$ is bijective. Thus $(x,z) \mapsto (1 - z A_\fa(x)^{-1})^{-1} \in L(H,H)$ is $\sC$   
  and hence also 
  \[
  (x,z) \mapsto (A_\fa(x)-z)^{-1}=A_\fa(x)^{-1}(1 - z A_\fa(x)^{-1})^{-1} \in L(H,H)
  \]
  is $\sC$ near $(x_0,z_0)$.
  This shows Claim~\ref{clnew1}.

  In what follows we assume that the parameter space is $E=\R$ and $t=x$. 
  
  \begin{claim}
  Assume that $t \mapsto \fa(t)$ is $C^{0,1}$.  
  A $C^{0,1}$-eigenvalue $t \mapsto \la_j(t)$ of $t \mapsto A_\fa(t)$ cannot accelerate to $\infty$ in finite time.
  \end{claim}
  
  Note that \eqref{faB} implies that $\fa(t)$ is locally uniformly sectorial. Thus we may assume without loss of 
  generality that $\fs(t) = \on{Re} a(t) \ge 1$ near $t_0$. Since $\fa(t)$ is closed, the inner product 
  $\<u \mid v\>_t := \<u \mid v\> + \fs(t)(u,v)$ makes $V$ to a Hilbert space $V_t:=(V,\|~\|_t)$ (see \cite[VI Thm.\ 1.11]{Kato76}). 
  The arguments in the proof of \cite[Claim (1)]{KMRp} show that all these norms $\|~\|_t$ are locally uniformly 
  equivalent.
  
  By reducing to $P_\fa(t)A_\fa(t)|_{P_\fa(t)(H)}$ (where $P_\fa(t) = -\frac{1}{2 \pi i} \int_\ga (A_\fa(x)-z)^{-1} dz$) 
  we have \eqref{abs} (with $t$ replaced by $t_0$) and hence, using \eqref{faB},  
  \begin{align*}
    \Big|&\frac{\la_j(t_k)-\la_j(t_0)}{t_k-t_0}\Big| |\< v_j(t_k) \mid w_j(t_0)\>|
    = \Big|\Big\<\frac{B(t_k)-B(t_0)}{t_k-t_0} G v_j(t_k) \mid G w_j(t_0)\Big\>\Big|\\
    &\le 
    \Big\|\frac{B(t_k)-B(t_0)}{t_k-t_0}\Big\|_{L(H,H)} \|G\|_{L(V_{t_k},H)}^2 \|v_j(t_k)\|_{V_{t_k}} \|w_j(t_0)\|_{V_{t_k}}\\
    &\le C   \sqrt{1 + \fs(t_k)(v_j(t_k))}  \cdot \sqrt{1 + \fs(t_0)(w_j(t_0))}
    =  C   \sqrt{1 + \on{Re} \la_j(t_k)}  \cdot \sqrt{1 + \on{Re} \la_j(t_0)},
  \end{align*}
  for a constant $C$, since all norms $\|~\|_t$ are locally uniformly equivalent 
  and since $\fs(t)(u,v) = \< S(t)^{1/2}u \mid S(t)^{1/2}v\>$ for $S(t)= A_{\fs(t)}\ge 1$,
  by the second representation theorem \cite[VI Thm.\ 2.23]{Kato76}. Since $t \mapsto \la_j(t)$ is Lipschitz, it follows that 
  \[
  |\on{Re} \la_j'(t)| \le |\la_j'(t)| \le C+C |\on{Re} \la_j(t)| \quad \text{ a.e.,}
  \]
  and Gronwall's lemma implies that $t \mapsto \on{Re} \la_j(t)$ cannot accelerate to $\infty$ in finite time.
  Since $A_\fa(t)$ is locally uniformly m-sectorial, $\la_j(t)$ lies in a sector $\{z \in \C : |\arg(z-\ze)| \le \th\}$, for 
  $0 \le \th < \pi/2$ and $\ze \in \R$, and the claim follows.  
\end{remark}

\section{The results are best possible} \label{optimal}

The condition on the order of contact in \thetag{\hyperref[A]{A}} cannot be dropped: 
This follows from the examples in \cite{KM03} and \cite{KMRp} for $C^\infty$ and for non-quasianalytic $C^{\{M\}}$.
From the latter one can also deduce a counterexample for non-quasianalytic $C^{(M)}$.
  
These examples together with Example~\ref{ex1} also show that
results of type \thetag{\hyperref[C]{C}}--\thetag{\hyperref[F]{F}} are hopeless for the eigenvectors.
Moreover,  
\thetag{\hyperref[B]{B}} is wrong without desingularization, by Example~\ref{ex1}.

Result \thetag{\hyperref[C]{C}} is optimal, since by Example~\ref{excont} the single eigenvalues cannot be chosen continuously in general.  
By the example in \cite{KM03}, in \thetag{\hyperref[E]{E}} and \thetag{\hyperref[F]{F}} the eigenvalues cannot be $C^{1,\be}$ 
for any $\be>0$, even if $t \mapsto A(t)$ is $C^\infty$. 
On the other hand, in our proof the assumption $C^{1,\al}$ in \thetag{\hyperref[E]{E}} (resp.\ $C^{2,\al}$ in \thetag{\hyperref[F]{F}}) cannot be weakened to $C^1$, 
by the ``resolvent example'' in \cite{KM03}, but we do not know whether there is a $C^1$ (resp.\ $C^2$) curve 
of unbounded normal operators with common domain and compact resolvent 
whose eigenvalues cannot be parameterized $C^1$ (resp.\ twice differentiably).

Example~\ref{ex3} and Example~\ref{ex4} show that the results are generally no longer true if $A$ is 
a family of merely diagonalizable matrices.

\begin{example}[The first partials of eigenvectors cannot be locally bounded] \label{ex1}
  The real analytic family of normal (even real symmetric) matrices
  \[
  A(x,y) = 
  \begin{pmatrix}
    x & y \\
    y & -x
  \end{pmatrix},
  \quad x,y \in \R,
  \]
  has the eigenvalues $\pm \sqrt{x^2+y^2}$. There cannot exist a parameterization of the eigenvectors
  of $A$ with locally bounded derivatives. Namely, if $\binom{u}{v}$ denotes an eigenvector with norm $1$ 
  for the eigenvalue $\sqrt{x^2+y^2}$, then the partial derivative $\binom{u_x}{v_x}$ (where it exists) must satisfy
  \[
  \begin{pmatrix}
    x-\sqrt{x^2+y^2} & y\\
    y & -x-\sqrt{x^2+y^2}
  \end{pmatrix}
  \binom{u_x}{v_x} 
  = 
  \begin{pmatrix}
    \frac{x}{\sqrt{x^2+y^2}}-1 & 0 \\
    0 & \frac{x}{\sqrt{x^2+y^2}}+1
  \end{pmatrix}
  \binom{u}{v}
  \]
  If $\binom{u_x}{v_x}$ were bounded near $0$, the left-hand side would converge to $0$ as $x,y \to 0$,
  whereas the right-hand side does not, a contradiction. 
\end{example}

\begin{example}[The single eigenvalues cannot be chosen continuously] \label{excont}
  The eigenvalues of the locally Lipschitz family of normal matrices
  \[
  A(x) = 
  \begin{pmatrix}
    0 & x \\
    |x| & 0
  \end{pmatrix},
  \quad x \in \C \cong \R^2
  \]
  do not admit a parameterization which is continuous in a neighborhood of $0$. 
\end{example}

\begin{example}[Mere diagonalizability does not guarantee $C^\om$ eigenvalues] \label{ex3}
  The eigenvalues of the real analytic curve of diagonalizable matrices
  \[
  A(x)= \begin{pmatrix}
    x & 0 & 0 \\
    0 & 0 & x^2 \\
    0 & x & 0
  \end{pmatrix},
  \quad x \in \R,
  \]
  are given by $x$, $\pm x^{3/2}$ if $x\ge 0$ and by $x$, $\pm i |x|^{3/2}$ if $x<0$.
\end{example}

\begin{example}[Mere diagonalizability does not guarantee $C^{0,1}$ eigenvalues] \label{ex4}
  See \cite[II Ex.\ 5.9]{Kato76}: 
  \[
  A(x) = 
  \begin{pmatrix}
    |x|^\al & |x|^\al-|x|^\be \big(2+\sin \frac{1}{|x|} \big) \\
    -|x|^\al & -|x|^\al
  \end{pmatrix}, 
  \quad x \in \R \setminus \{0\}, \quad A(0)=0,
  \]
  forms a $C^1$-curve of diagonalizable matrices if $\al>1$ and $\be>2$.
  The eigenvalues of $A$ are given by 
  \[
  \la_\pm(x) = \pm |x|^{\frac{\al+\be}{2}} \Big(2+\sin \frac{1}{|x|} \Big)^{\frac{1}{2}}, \quad x \in \R \setminus \{0\}, \quad \la_\pm(0)=0.
  \] 
  The derivatives $\la_\pm'$ exist everywhere, but they are discontinuous at $0$ if $\al +\be \le 4$ and even unbounded near $0$ if $\al +\be < 4$.   
\end{example}

\subsection*{Acknowledgements} I am grateful to Adam Parusi\'nski for helpful discussions and to the 
anonymous referee for pointing to m-sectorial operators.

\def\cprime{$'$}
\providecommand{\bysame}{\leavevmode\hbox to3em{\hrulefill}\thinspace}
\providecommand{\MR}{\relax\ifhmode\unskip\space\fi MR }
\providecommand{\MRhref}[2]{%
  \href{http://www.ams.org/mathscinet-getitem?mr=#1}{#2}
}
\providecommand{\href}[2]{#2}

\end{document}